\documentclass[12pt,a4paper]{article}

\usepackage{theorem,enumerate}
\usepackage{amsmath,latexsym,amssymb,amsfonts}
\usepackage{eucal}
\usepackage{color}
\usepackage{comment}

\newcounter{Scounter}
\theorembodyfont{\normalfont\slshape}
\newtheorem{thm}{Theorem}[section]

\newtheorem{Thm}{Theorem}

\newtheorem{lem}[thm]{Lemma}

\newtheorem{Q}[thm]{Question} 

\newtheorem{claim}{Claim}[section]

\newtheorem{con}{Conjecture}

\numberwithin{equation}{section}
\newtheorem{remark}{\normalfont\itshape Remark}

\newcommand{\proof}{\medbreak\noindent\textit{Proof.}\quad}

\newcommand{\qed}{{$\quad\square$\vs{3.6}}}

\newcommand{\vs}[1]{\vspace*{#1 mm}}

\def\C{{ \mathcal{C}}}

\def\F{{ \mathcal{F}}}
\def\G{{ \mathcal{G}}}
\def\H{{ \mathcal{H}}}

\def\P{{ \mathcal{P}}}

\bfseries\normalfont

\title{Difference of forbidden pairs containing a claw}

\author{
Guantao Chen$^{1}$\footnote{This research partially supported by an NSA grant}
\and
Michitaka Furuya$^{2}$\footnote{Email address: \texttt{michitaka.furuya@gmail.com}, This research was partially supported by JSPS KAKENHI Grant number 26800086}\and
Songling Shan$^{1}$ \and
Shoichi Tsuchiya$^{3}$\footnote{Email address: \texttt{s.tsuchiya@isc.senshu-u.ac.jp}}\qquad
Ping Yang$^{1}$
\\
\small
$^{1}$\small\textsl{Dept. of Math \& Stat,
Georgia State University, Atlanta}\\ 
$^{2}$\small\textsl{Dept. of Mathematical Information Science, 
Tokyo University of Science, Tokyo }\\ 
$^{3}$\small\textsl{School of Network and Information, Senshu University, Kanagawa}\\ 
}

\begin{document}

\maketitle

\begin{abstract}
When we study forbidden subgraph conditions guaranteeing graphs to have some properties, a claw (or $K_{1,3}$) frequently appears as one of forbidden subgraphs.
Recently, Furuya and Tsuchiya compared two classes generated by different forbidden pairs containing a claw, and characterized one of such classes.
In this paper, we give such characterization for three new classes.
Furthermore, we give applications of our characterizations to some forbidden subgraph problems.
\end{abstract}

\noindent
{\it Key words and phrases.}
forbidden subgraph, Hamiltonian cycle, Halin graph.

\noindent
{\it AMS 2010 Mathematics Subject Classification.}
05C75.

\section{Introduction}\label{sec1}

Let $\G_{1}$ and $\G_{2}$ be two families of graphs, and let $P$ be a certain property for graphs.
We assume that every member of $\G_{2}$ satisfies $P$, and consider the problem whether members of $\G_{1}$ satisfy $P$ or not.
If we suppose $\G_{1}\subseteq \G_{2}$, then every member of $\G_{1}$ satisfies $P$.
Now, we suppose a weaker condition than $\G_{1}\subseteq \G_{2}$.

We first suppose that the family $\G_{1}-\G_{2}$ is finite.
Then every member of $\G_{1}$ satisfies $P$ with finite exceptions.
Since we can check whether finite members of $\G_{1}$ satisfy $P$ or not in finite time, we can regard the desired problem as solved.

We next suppose that the members of $\G_{1}-\G_{2}$ is characterized (not necessary finite).
Then each member of $\G_{1}$ either satisfies $P$ or is characterized.
If the characterization has a simple structure, then we may be able to check whether such graphs satisfy $P$ or not.
Thus, in this case, it might be possible to solve the desired problem.

In this paper, we try to apply the above strategy for the forbidden subgraph problems.

\subsection{Definition and preliminary}\label{subsec1.1}

For a family $\F $ of connected graphs, a graph $G$ is said to be {\it $\F $-free} if $G$ contains no member of $\F $ as an induced subgraph.
We also say that the members of $\F $ are {\it forbidden subgraphs}.
If $G$ is $\{F\}$-free, then $G$ is simply said to be {\it $F$-free}.
A family $\F$ of forbidden subgraphs is called a {\it forbidden pair} if $|\F|=2$.

\begin{figure}[htb]
\begin{center}
{\unitlength 0.1in
\begin{picture}( 43.0000,  8.2500)(  3.5000,-10.3500)
%
\special{sh 1.000}%
\special{ia 600 460 50 50  0.0000000  6.2831853}%
\special{pn 8}%
\special{ar 600 460 50 50  0.0000000  6.2831853}%
%
\special{sh 1.000}%
\special{ia 400 800 50 50  0.0000000  6.2831853}%
\special{pn 8}%
\special{ar 400 800 50 50  0.0000000  6.2831853}%
%
\special{sh 1.000}%
\special{ia 800 800 50 50  0.0000000  6.2831853}%
\special{pn 8}%
\special{ar 800 800 50 50  0.0000000  6.2831853}%
%
\special{pn 8}%
\special{pa 600 460}%
\special{pa 400 800}%
\special{fp}%
\special{pa 400 800}%
\special{pa 800 800}%
\special{fp}%
\special{pa 800 800}%
\special{pa 600 460}%
\special{fp}%
%
\special{sh 1.000}%
\special{ia 990 800 50 50  0.0000000  6.2831853}%
\special{pn 8}%
\special{ar 990 800 50 50  0.0000000  6.2831853}%
%
\special{sh 1.000}%
\special{ia 1390 800 50 50  0.0000000  6.2831853}%
\special{pn 8}%
\special{ar 1390 800 50 50  0.0000000  6.2831853}%
%
\special{pn 4}%
\special{sh 1}%
\special{ar 1190 800 8 8 0  6.28318530717959E+0000}%
\special{sh 1}%
\special{ar 1240 800 8 8 0  6.28318530717959E+0000}%
\special{sh 1}%
\special{ar 1140 800 8 8 0  6.28318530717959E+0000}%
\special{sh 1}%
\special{ar 1140 800 8 8 0  6.28318530717959E+0000}%
%
\special{pn 8}%
\special{pa 800 800}%
\special{pa 1090 800}%
\special{fp}%
\special{pa 1300 800}%
\special{pa 1300 800}%
\special{fp}%
\special{pa 1300 800}%
\special{pa 1390 800}%
\special{fp}%
%
\special{pn 8}%
\special{pa 1006 710}%
\special{pa 1006 610}%
\special{fp}%
\special{pa 1006 660}%
\special{pa 1146 660}%
\special{fp}%
\special{pa 1246 660}%
\special{pa 1406 660}%
\special{fp}%
%
\special{pn 8}%
\special{pa 1406 710}%
\special{pa 1406 610}%
\special{fp}%
\put(12.0500,-5.0000){\makebox(0,0){$n$ vertices}}%
%
\special{pn 8}%
\special{pa 1406 660}%
\special{pa 1006 660}%
\special{fp}%
%
\special{sh 1.000}%
\special{ia 2200 460 50 50  0.0000000  6.2831853}%
\special{pn 8}%
\special{ar 2200 460 50 50  0.0000000  6.2831853}%
%
\special{sh 1.000}%
\special{ia 2000 800 50 50  0.0000000  6.2831853}%
\special{pn 8}%
\special{ar 2000 800 50 50  0.0000000  6.2831853}%
%
\special{sh 1.000}%
\special{ia 2400 800 50 50  0.0000000  6.2831853}%
\special{pn 8}%
\special{ar 2400 800 50 50  0.0000000  6.2831853}%
%
\special{pn 8}%
\special{pa 2200 460}%
\special{pa 2000 800}%
\special{fp}%
\special{pa 2000 800}%
\special{pa 2400 800}%
\special{fp}%
\special{pa 2400 800}%
\special{pa 2200 460}%
\special{fp}%
%
\special{sh 1.000}%
\special{ia 2590 800 50 50  0.0000000  6.2831853}%
\special{pn 8}%
\special{ar 2590 800 50 50  0.0000000  6.2831853}%
%
\special{sh 1.000}%
\special{ia 2990 800 50 50  0.0000000  6.2831853}%
\special{pn 8}%
\special{ar 2990 800 50 50  0.0000000  6.2831853}%
%
\special{pn 4}%
\special{sh 1}%
\special{ar 2790 800 8 8 0  6.28318530717959E+0000}%
\special{sh 1}%
\special{ar 2840 800 8 8 0  6.28318530717959E+0000}%
\special{sh 1}%
\special{ar 2740 800 8 8 0  6.28318530717959E+0000}%
\special{sh 1}%
\special{ar 2740 800 8 8 0  6.28318530717959E+0000}%
%
\special{pn 8}%
\special{pa 2400 800}%
\special{pa 2690 800}%
\special{fp}%
\special{pa 2900 800}%
\special{pa 2900 800}%
\special{fp}%
\special{pa 2900 800}%
\special{pa 2990 800}%
\special{fp}%
\put(8.0000,-11.0000){\makebox(0,0){$Z_{n}$}}%
\put(24.0000,-11.0000){\makebox(0,0){$B_{1,n}$}}%
%
\special{sh 1.000}%
\special{ia 2200 260 50 50  0.0000000  6.2831853}%
\special{pn 8}%
\special{ar 2200 260 50 50  0.0000000  6.2831853}%
%
\special{pn 8}%
\special{pa 2200 260}%
\special{pa 2200 460}%
\special{fp}%
%
\special{sh 1.000}%
\special{ia 4200 460 50 50  0.0000000  6.2831853}%
\special{pn 8}%
\special{ar 4200 460 50 50  0.0000000  6.2831853}%
%
\special{sh 1.000}%
\special{ia 4000 800 50 50  0.0000000  6.2831853}%
\special{pn 8}%
\special{ar 4000 800 50 50  0.0000000  6.2831853}%
%
\special{sh 1.000}%
\special{ia 4400 800 50 50  0.0000000  6.2831853}%
\special{pn 8}%
\special{ar 4400 800 50 50  0.0000000  6.2831853}%
%
\special{pn 8}%
\special{pa 4200 460}%
\special{pa 4000 800}%
\special{fp}%
\special{pa 4000 800}%
\special{pa 4400 800}%
\special{fp}%
\special{pa 4400 800}%
\special{pa 4200 460}%
\special{fp}%
%
\special{sh 1.000}%
\special{ia 4200 260 50 50  0.0000000  6.2831853}%
\special{pn 8}%
\special{ar 4200 260 50 50  0.0000000  6.2831853}%
%
\special{sh 1.000}%
\special{ia 4600 800 50 50  0.0000000  6.2831853}%
\special{pn 8}%
\special{ar 4600 800 50 50  0.0000000  6.2831853}%
%
\special{sh 1.000}%
\special{ia 3800 800 50 50  0.0000000  6.2831853}%
\special{pn 8}%
\special{ar 3800 800 50 50  0.0000000  6.2831853}%
%
\special{pn 8}%
\special{pa 3800 800}%
\special{pa 4000 800}%
\special{fp}%
\special{pa 4400 800}%
\special{pa 4600 800}%
\special{fp}%
%
\special{pn 8}%
\special{pa 4200 260}%
\special{pa 4200 460}%
\special{fp}%
\put(42.0000,-11.0000){\makebox(0,0){$N$}}%
%
\special{pn 8}%
\special{pa 2606 706}%
\special{pa 2606 606}%
\special{fp}%
\special{pa 2606 656}%
\special{pa 2746 656}%
\special{fp}%
\special{pa 2846 656}%
\special{pa 3006 656}%
\special{fp}%
%
\special{pn 8}%
\special{pa 3006 706}%
\special{pa 3006 606}%
\special{fp}%
\put(28.0500,-4.9500){\makebox(0,0){$n$ vertices}}%
%
\special{pn 8}%
\special{pa 3006 656}%
\special{pa 2606 656}%
\special{fp}%
\end{picture}}%
\caption{Graphs $Z_{n}$, $B_{1,n}$ and $N$.}
\label{fig_N}
\end{center}
\end{figure}
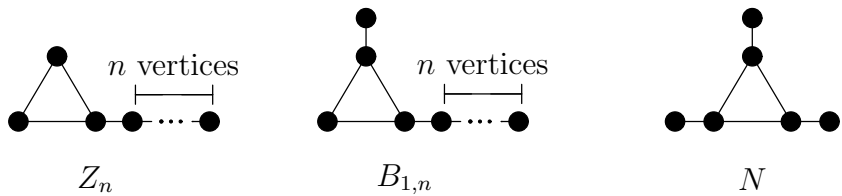

Let $K_{1,3}$ denote the star with three leaves. Let $K_{n}$ and $P_{n}$ denote the complete graph and the path of order $n$, respectively.
For nonnegative integers $k$, $l$ and $m$, let $N_{k,l,m}$ be a graph obtained from $K_3$ and vertex disjoint three paths $P_{k+1}$, $P_{l+1}$, $P_{m+1}$ by identifying one end-vertex of the paths and distinct three vertices of the $K_3$.
Commonly, $N_{k,0,0}$ (resp., $N_{k,l,0}$) is usually denoted by $Z_{k}$ (resp., $B_{k,l}$), 
and $N_{1,1,1}$ is usually denoted by $N$ (see Figure~\ref{fig_N}).

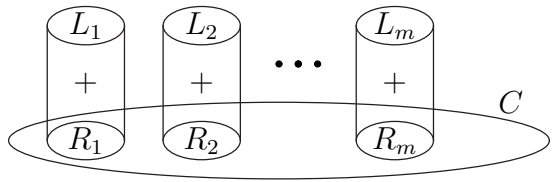
\begin{figure}[htb]
\begin{center}
{\unitlength 0.1in
\begin{picture}( 28.0000,  9.0000)(  4.0000,-12.0000)
%
\special{pn 8}%
\special{ar 800 400 200 100  0.0000000  6.2831853}%
%
\special{pn 8}%
\special{ar 1400 400 200 100  0.0000000  6.2831853}%
%
\special{pn 8}%
\special{ar 2400 400 200 100  0.0000000  6.2831853}%
%
\special{pn 8}%
\special{ar 1800 1000 1400 200  0.0000000  6.2831853}%
%
\special{pn 8}%
\special{ar 800 1000 200 100  0.0000000  6.2831853}%
%
\special{pn 8}%
\special{ar 1400 1000 200 100  0.0000000  6.2831853}%
%
\special{pn 8}%
\special{ar 2400 1000 200 100  0.0000000  6.2831853}%
%
\special{pn 8}%
\special{pa 600 1000}%
\special{pa 600 400}%
\special{fp}%
\special{pa 1000 1000}%
\special{pa 1000 400}%
\special{fp}%
\special{pa 1200 400}%
\special{pa 1200 1000}%
\special{fp}%
\special{pa 1600 1000}%
\special{pa 1600 400}%
\special{fp}%
\special{pa 2200 400}%
\special{pa 2200 1000}%
\special{fp}%
\special{pa 2600 1000}%
\special{pa 2600 400}%
\special{fp}%
%
\special{pn 4}%
\special{sh 1}%
\special{ar 1800 600 16 16 0  6.28318530717959E+0000}%
\special{sh 1}%
\special{ar 2000 600 16 16 0  6.28318530717959E+0000}%
\special{sh 1}%
\special{ar 1900 600 16 16 0  6.28318530717959E+0000}%
\special{sh 1}%
\special{ar 1900 600 16 16 0  6.28318530717959E+0000}%
\put(8.0000,-4.0000){\makebox(0,0){$L_{1}$}}%
\put(14.0000,-4.0000){\makebox(0,0){$L_{2}$}}%
\put(24.0000,-4.0000){\makebox(0,0){$L_{m}$}}%
\put(8.0000,-10.0000){\makebox(0,0){$R_{1}$}}%
\put(14.0000,-10.0000){\makebox(0,0){$R_{2}$}}%
\put(24.0000,-10.0000){\makebox(0,0){$R_{m}$}}%
\put(30.0000,-8.0000){\makebox(0,0){$C$}}%
\put(8.0000,-7.0000){\makebox(0,0){$+$}}%
\put(14.0000,-7.0000){\makebox(0,0){$+$}}%
\put(24.0000,-7.0000){\makebox(0,0){$+$}}%
\end{picture}}%

\caption{Generalized comb}
\label{f0}
\end{center}
\end{figure}

As we mentioned, our main aim is to characterize connected $\F_{1}$-free but not $\F_{2}$-free graphs for two families $\F_{1}$ and $\F_{2}$ of forbidden subgraphs.
Such study derives from \cite{O}.
In~\cite{O}, Olariu considered the case where $\F_{1}=\{Z_{1}\}$ and $\F_{2}=\{K_{3}\}$, and showed that every connected $Z_{1}$-free but not $K_{3}$-free graph is a complete multipartite graph with at least three partite sets.
The result is useful when we investigate the class of $Z_{1}$-free graphs (for example, the characterization was used for research of perfect $Z_{1}$-free graphs in \cite{O}).
Recently, Furuya and Tsuchiya~\cite{FT} focused on forbidden pairs for the existence of a Hamiltonian cycle, and studied characterization similar to Olariu's result.
A graph $H$ is a {\it generalized comb} if $H$ is obtained as follows (see Figure~\ref{f0}):
Let $m\geq 3$ be an integer.
Let $L_{i}~(1\leq i\leq m)$ and $C$ be vertex-disjoint non-empty cliques with $|C|\geq m$, and let $R_{i}~(1\leq i\leq m)$ be disjoint non-empty subcliques of $C$.
We define the graph $H$ on $\Big(\bigcup _{1\leq i\leq m}L_{i}\Big) \cup C$ such that every vertex in $L_{i}$ is joined to all vertices in $R_{i}$ for each $i~(1\leq i\leq m)$.
In this context, $L_{i}$ is called a {\it leaf-clique} and $R_{i}$ is called the {\it root} of $L_{i}$.
The following theorem was proved in \cite{FT}.

\begin{Thm}[Furuya and Tsuchiya~\cite{FT}]
\label{ThmA}
Let $G$ be a connected $\{K_{1,3},B_{1,2}\}$-free graph.
Then $G$ is not $N$-free if and only if $G$ is a generalized comb.
\end{Thm}

In other words, they solved a characterization problem for $\F_{1}=\{K_{1,3},B_{1,2}\}$ and $\F_{2}=\{K_{1,3},N\}$.

Our notation and terminology are standard, and mostly taken from \cite{D}.
In particular, we shall use the following terminology.
Let $G$ be a graph.
For $v\in V(G)$, we let $N_{G}(v)$ denote the {\it neighborhood} of $v$ in $G$.
For a set $U$, we let $G[U]$ denote the subgraph of $G$ induced by $U\cap V(G)$.

\subsection{Main results}\label{subsec1.2}

In this paper, we investigate graphs generated by different families of forbidden subgraphs, and characterize the following classes:
\begin{enumerate}[(F1)]
\item
connected $\{K_{1,3},Z_{2}\}$-free but not $B_{1,1}$-free graphs,
\item
connected $\{K_{1,3},B_{1,1}\}$-free but not $P_{5}$-free graphs, and
\item
connected $\{K_{1,3},B_{1,2}\}$-free but not $P_{6}$-free graphs.
\end{enumerate}

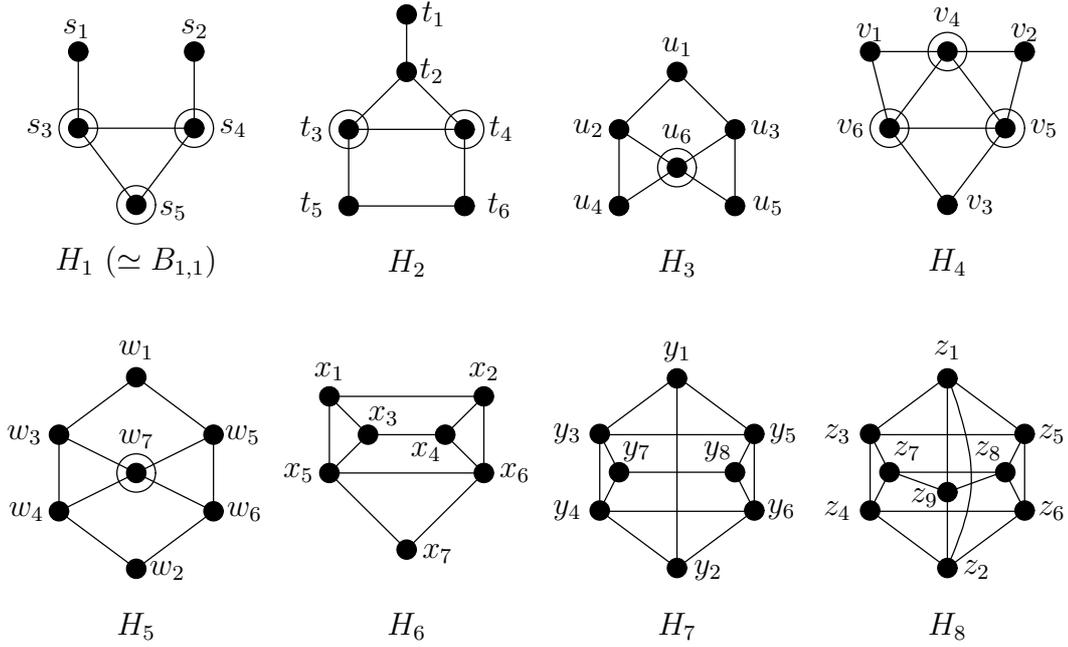
\begin{figure}[htb]
\begin{center}
\unitlength 0.1in
\begin{picture}( 54.4000, 31.9500)( -0.9000,-33.3500)
%
\special{sh 1.000}%
\special{ia 400 800 50 50  0.0000000  6.2831853}%
\special{pn 8}%
\special{pn 8}%
\special{ar 400 800 50 50  0.0000000  6.2831853}%
%
\special{sh 1.000}%
\special{ia 1000 800 50 50  0.0000000  6.2831853}%
\special{pn 8}%
\special{pn 8}%
\special{ar 1000 800 50 50  0.0000000  6.2831853}%
%
\special{sh 1.000}%
\special{ia 700 1200 50 50  0.0000000  6.2831853}%
\special{pn 8}%
\special{pn 8}%
\special{ar 700 1200 50 50  0.0000000  6.2831853}%
%
\special{sh 1.000}%
\special{ia 400 400 50 50  0.0000000  6.2831853}%
\special{pn 8}%
\special{pn 8}%
\special{ar 400 400 50 50  0.0000000  6.2831853}%
%
\special{sh 1.000}%
\special{ia 1000 400 50 50  0.0000000  6.2831853}%
\special{pn 8}%
\special{pn 8}%
\special{ar 1000 400 50 50  0.0000000  6.2831853}%
%
\special{pn 8}%
\special{pa 1000 400}%
\special{pa 1000 800}%
\special{fp}%
\special{pa 1000 800}%
\special{pa 400 800}%
\special{fp}%
\special{pa 400 800}%
\special{pa 400 400}%
\special{fp}%
\special{pa 400 800}%
\special{pa 700 1200}%
\special{fp}%
\special{pa 700 1200}%
\special{pa 1000 800}%
\special{fp}%
%
\special{pn 8}%
\special{ar 400 800 100 100  0.0000000  6.2831853}%
%
\special{pn 8}%
\special{ar 1000 800 100 100  0.0000000  6.2831853}%
%
\special{pn 8}%
\special{ar 700 1200 100 100  0.0000000  6.2831853}%
%
\special{sh 1.000}%
\special{ia 4600 800 50 50  0.0000000  6.2831853}%
\special{pn 8}%
\special{pn 8}%
\special{ar 4600 800 50 50  0.0000000  6.2831853}%
%
\special{sh 1.000}%
\special{ia 5200 800 50 50  0.0000000  6.2831853}%
\special{pn 8}%
\special{pn 8}%
\special{ar 5200 800 50 50  0.0000000  6.2831853}%
%
\special{sh 1.000}%
\special{ia 4900 400 50 50  0.0000000  6.2831853}%
\special{pn 8}%
\special{pn 8}%
\special{ar 4900 400 50 50  0.0000000  6.2831853}%
%
\special{sh 1.000}%
\special{ia 5300 400 50 50  0.0000000  6.2831853}%
\special{pn 8}%
\special{pn 8}%
\special{ar 5300 400 50 50  0.0000000  6.2831853}%
%
\special{sh 1.000}%
\special{ia 4500 400 50 50  0.0000000  6.2831853}%
\special{pn 8}%
\special{pn 8}%
\special{ar 4500 400 50 50  0.0000000  6.2831853}%
%
\special{sh 1.000}%
\special{ia 4900 1200 50 50  0.0000000  6.2831853}%
\special{pn 8}%
\special{pn 8}%
\special{ar 4900 1200 50 50  0.0000000  6.2831853}%
%
\special{pn 8}%
\special{pa 4900 1200}%
\special{pa 4600 800}%
\special{fp}%
\special{pa 4600 800}%
\special{pa 5200 800}%
\special{fp}%
\special{pa 5200 800}%
\special{pa 5300 400}%
\special{fp}%
\special{pa 5300 400}%
\special{pa 4500 400}%
\special{fp}%
\special{pa 4500 400}%
\special{pa 4600 800}%
\special{fp}%
\special{pa 4600 800}%
\special{pa 4900 400}%
\special{fp}%
\special{pa 4900 400}%
\special{pa 5200 800}%
\special{fp}%
\special{pa 5200 800}%
\special{pa 4900 1200}%
\special{fp}%
%
\special{pn 8}%
\special{ar 4900 400 100 100  0.0000000  6.2831853}%
%
\special{pn 8}%
\special{ar 4600 800 100 100  0.0000000  6.2831853}%
%
\special{pn 8}%
\special{ar 5200 800 100 100  0.0000000  6.2831853}%
%
\special{sh 1.000}%
\special{ia 1800 1206 50 50  0.0000000  6.2831853}%
\special{pn 8}%
\special{pn 8}%
\special{ar 1800 1206 50 50  0.0000000  6.2831853}%
%
\special{sh 1.000}%
\special{ia 2400 1206 50 50  0.0000000  6.2831853}%
\special{pn 8}%
\special{pn 8}%
\special{ar 2400 1206 50 50  0.0000000  6.2831853}%
%
\special{sh 1.000}%
\special{ia 2400 806 50 50  0.0000000  6.2831853}%
\special{pn 8}%
\special{pn 8}%
\special{ar 2400 806 50 50  0.0000000  6.2831853}%
%
\special{sh 1.000}%
\special{ia 1800 806 50 50  0.0000000  6.2831853}%
\special{pn 8}%
\special{pn 8}%
\special{ar 1800 806 50 50  0.0000000  6.2831853}%
%
\special{sh 1.000}%
\special{ia 2100 506 50 50  0.0000000  6.2831853}%
\special{pn 8}%
\special{pn 8}%
\special{ar 2100 506 50 50  0.0000000  6.2831853}%
%
\special{sh 1.000}%
\special{ia 2100 206 50 50  0.0000000  6.2831853}%
\special{pn 8}%
\special{pn 8}%
\special{ar 2100 206 50 50  0.0000000  6.2831853}%
%
\special{pn 8}%
\special{pa 2100 206}%
\special{pa 2100 506}%
\special{fp}%
\special{pa 2100 506}%
\special{pa 1800 806}%
\special{fp}%
\special{pa 1800 806}%
\special{pa 2400 806}%
\special{fp}%
\special{pa 2400 806}%
\special{pa 2400 1206}%
\special{fp}%
\special{pa 2400 1206}%
\special{pa 1800 1206}%
\special{fp}%
\special{pa 1800 1206}%
\special{pa 1800 806}%
\special{fp}%
\special{pa 2400 806}%
\special{pa 2100 506}%
\special{fp}%
%
\special{pn 8}%
\special{ar 1800 806 100 100  0.0000000  6.2831853}%
%
\special{pn 8}%
\special{ar 2400 806 100 100  0.0000000  6.2831853}%
%
\special{sh 1.000}%
\special{ia 3200 806 50 50  0.0000000  6.2831853}%
\special{pn 8}%
\special{pn 8}%
\special{ar 3200 806 50 50  0.0000000  6.2831853}%
%
\special{sh 1.000}%
\special{ia 3200 1206 50 50  0.0000000  6.2831853}%
\special{pn 8}%
\special{pn 8}%
\special{ar 3200 1206 50 50  0.0000000  6.2831853}%
%
\special{sh 1.000}%
\special{ia 3800 806 50 50  0.0000000  6.2831853}%
\special{pn 8}%
\special{pn 8}%
\special{ar 3800 806 50 50  0.0000000  6.2831853}%
%
\special{sh 1.000}%
\special{ia 3800 1206 50 50  0.0000000  6.2831853}%
\special{pn 8}%
\special{pn 8}%
\special{ar 3800 1206 50 50  0.0000000  6.2831853}%
%
\special{sh 1.000}%
\special{ia 3500 1006 50 50  0.0000000  6.2831853}%
\special{pn 8}%
\special{pn 8}%
\special{ar 3500 1006 50 50  0.0000000  6.2831853}%
%
\special{sh 1.000}%
\special{ia 3500 506 50 50  0.0000000  6.2831853}%
\special{pn 8}%
\special{pn 8}%
\special{ar 3500 506 50 50  0.0000000  6.2831853}%
%
\special{pn 8}%
\special{pa 3500 506}%
\special{pa 3200 806}%
\special{fp}%
\special{pa 3200 806}%
\special{pa 3200 1206}%
\special{fp}%
\special{pa 3200 1206}%
\special{pa 3500 1006}%
\special{fp}%
\special{pa 3500 1006}%
\special{pa 3200 806}%
\special{fp}%
\special{pa 3500 1006}%
\special{pa 3800 1206}%
\special{fp}%
\special{pa 3800 1206}%
\special{pa 3800 806}%
\special{fp}%
\special{pa 3800 806}%
\special{pa 3500 506}%
\special{fp}%
\special{pa 3500 1006}%
\special{pa 3800 806}%
\special{fp}%
%
\special{pn 8}%
\special{ar 3500 1006 100 100  0.0000000  6.2831853}%
%
\special{sh 1.000}%
\special{ia 1700 2200 50 50  0.0000000  6.2831853}%
\special{pn 8}%
\special{pn 8}%
\special{ar 1700 2200 50 50  0.0000000  6.2831853}%
%
\special{sh 1.000}%
\special{ia 1700 2600 50 50  0.0000000  6.2831853}%
\special{pn 8}%
\special{pn 8}%
\special{ar 1700 2600 50 50  0.0000000  6.2831853}%
%
\special{sh 1.000}%
\special{ia 1900 2400 50 50  0.0000000  6.2831853}%
\special{pn 8}%
\special{pn 8}%
\special{ar 1900 2400 50 50  0.0000000  6.2831853}%
%
\special{sh 1.000}%
\special{ia 2300 2400 50 50  0.0000000  6.2831853}%
\special{pn 8}%
\special{pn 8}%
\special{ar 2300 2400 50 50  0.0000000  6.2831853}%
%
\special{sh 1.000}%
\special{ia 2500 2200 50 50  0.0000000  6.2831853}%
\special{pn 8}%
\special{pn 8}%
\special{ar 2500 2200 50 50  0.0000000  6.2831853}%
%
\special{sh 1.000}%
\special{ia 2500 2600 50 50  0.0000000  6.2831853}%
\special{pn 8}%
\special{pn 8}%
\special{ar 2500 2600 50 50  0.0000000  6.2831853}%
%
\special{sh 1.000}%
\special{ia 2100 3000 50 50  0.0000000  6.2831853}%
\special{pn 8}%
\special{pn 8}%
\special{ar 2100 3000 50 50  0.0000000  6.2831853}%
%
\special{pn 8}%
\special{pa 2100 3000}%
\special{pa 1700 2600}%
\special{fp}%
\special{pa 1700 2600}%
\special{pa 2500 2600}%
\special{fp}%
\special{pa 2500 2600}%
\special{pa 2100 3000}%
\special{fp}%
\special{pa 1900 2400}%
\special{pa 1700 2600}%
\special{fp}%
\special{pa 1700 2600}%
\special{pa 1700 2200}%
\special{fp}%
\special{pa 1700 2200}%
\special{pa 1900 2400}%
\special{fp}%
\special{pa 1900 2400}%
\special{pa 2300 2400}%
\special{fp}%
\special{pa 2300 2400}%
\special{pa 2500 2600}%
\special{fp}%
\special{pa 2500 2600}%
\special{pa 2500 2200}%
\special{fp}%
\special{pa 2500 2200}%
\special{pa 2300 2400}%
\special{fp}%
\special{pa 2500 2200}%
\special{pa 1700 2200}%
\special{fp}%
%
\special{sh 1.000}%
\special{ia 700 2100 50 50  0.0000000  6.2831853}%
\special{pn 8}%
\special{pn 8}%
\special{ar 700 2100 50 50  0.0000000  6.2831853}%
%
\special{sh 1.000}%
\special{ia 300 2400 50 50  0.0000000  6.2831853}%
\special{pn 8}%
\special{pn 8}%
\special{ar 300 2400 50 50  0.0000000  6.2831853}%
%
\special{sh 1.000}%
\special{ia 1100 2400 50 50  0.0000000  6.2831853}%
\special{pn 8}%
\special{pn 8}%
\special{ar 1100 2400 50 50  0.0000000  6.2831853}%
%
\special{sh 1.000}%
\special{ia 1100 2800 50 50  0.0000000  6.2831853}%
\special{pn 8}%
\special{pn 8}%
\special{ar 1100 2800 50 50  0.0000000  6.2831853}%
%
\special{sh 1.000}%
\special{ia 300 2800 50 50  0.0000000  6.2831853}%
\special{pn 8}%
\special{pn 8}%
\special{ar 300 2800 50 50  0.0000000  6.2831853}%
%
\special{sh 1.000}%
\special{ia 700 3100 50 50  0.0000000  6.2831853}%
\special{pn 8}%
\special{pn 8}%
\special{ar 700 3100 50 50  0.0000000  6.2831853}%
%
\special{pn 8}%
\special{pa 700 3100}%
\special{pa 300 2800}%
\special{fp}%
\special{pa 300 2800}%
\special{pa 300 2400}%
\special{fp}%
\special{pa 300 2400}%
\special{pa 700 2100}%
\special{fp}%
\special{pa 700 2100}%
\special{pa 1100 2400}%
\special{fp}%
\special{pa 1100 2400}%
\special{pa 1100 2800}%
\special{fp}%
\special{pa 1100 2800}%
\special{pa 700 3100}%
\special{fp}%
%
\special{sh 1.000}%
\special{ia 700 2600 50 50  0.0000000  6.2831853}%
\special{pn 8}%
\special{pn 8}%
\special{ar 700 2600 50 50  0.0000000  6.2831853}%
%
\special{pn 8}%
\special{pa 300 2400}%
\special{pa 1100 2800}%
\special{fp}%
\special{pa 1100 2400}%
\special{pa 300 2800}%
\special{fp}%
%
\special{pn 8}%
\special{ar 700 2600 100 100  0.0000000  6.2831853}%
\put(3.9000,-2.8000){\makebox(0,0){$s_{1}$}}%
\put(10.0000,-2.8000){\makebox(0,0){$s_{2}$}}%
\put(2.0000,-8.0000){\makebox(0,0){$s_{3}$}}%
\put(12.0000,-8.0000){\makebox(0,0){$s_{4}$}}%
\put(45.0000,-3.0000){\makebox(0,0){$v_{1}$}}%
\put(8.9000,-12.3000){\makebox(0,0){$s_{5}$}}%
\put(53.0000,-3.0000){\makebox(0,0){$v_{2}$}}%
\put(50.7000,-12.0000){\makebox(0,0){$v_{3}$}}%
\put(49.0000,-2.2000){\makebox(0,0){$v_{4}$}}%
\put(54.0000,-8.0000){\makebox(0,0){$v_{5}$}}%
\put(44.0000,-8.0000){\makebox(0,0){$v_{6}$}}%
\put(22.4000,-2.0500){\makebox(0,0){$t_{1}$}}%
\put(22.3000,-5.0000){\makebox(0,0){$t_{2}$}}%
\put(16.1000,-8.0500){\makebox(0,0){$t_{3}$}}%
\put(25.9000,-8.0500){\makebox(0,0){$t_{4}$}}%
\put(25.8000,-12.0500){\makebox(0,0){$t_{6}$}}%
\put(16.1000,-12.0000){\makebox(0,0){$t_{5}$}}%
\put(35.0000,-3.7500){\makebox(0,0){$u_{1}$}}%
\put(30.4000,-8.0000){\makebox(0,0){$u_{2}$}}%
\put(39.6000,-8.0500){\makebox(0,0){$u_{3}$}}%
\put(39.7000,-12.0500){\makebox(0,0){$u_{5}$}}%
\put(30.4000,-12.0000){\makebox(0,0){$u_{4}$}}%
\put(35.0000,-8.3000){\makebox(0,0){$u_{6}$}}%
\put(7.0000,-19.6000){\makebox(0,0){$w_{1}$}}%
\put(8.6000,-31.1000){\makebox(0,0){$w_{2}$}}%
\put(1.2000,-24.0000){\makebox(0,0){$w_{3}$}}%
\put(1.3000,-28.0000){\makebox(0,0){$w_{4}$}}%
\put(12.5000,-24.0000){\makebox(0,0){$w_{5}$}}%
\put(12.6000,-28.0000){\makebox(0,0){$w_{6}$}}%
\put(7.0000,-24.2000){\makebox(0,0){$w_{7}$}}%
\put(17.0000,-20.7000){\makebox(0,0){$x_{1}$}}%
\put(25.0000,-20.7000){\makebox(0,0){$x_{2}$}}%
\put(19.8000,-23.0000){\makebox(0,0){$x_{3}$}}%
\put(22.0000,-24.9000){\makebox(0,0){$x_{4}$}}%
\put(15.4000,-26.0000){\makebox(0,0){$x_{5}$}}%
\put(26.6000,-26.0000){\makebox(0,0){$x_{6}$}}%
\put(22.6000,-30.2000){\makebox(0,0){$x_{7}$}}%
\put(7.0000,-15.0000){\makebox(0,0){$H_{1}~(\simeq B_{1,1})$}}%
\put(49.0000,-15.0000){\makebox(0,0){$H_{4}$}}%
\put(21.0000,-15.0500){\makebox(0,0){$H_{2}$}}%
\put(35.0000,-15.0500){\makebox(0,0){$H_{3}$}}%
\put(7.0000,-34.0000){\makebox(0,0){$H_{5}$}}%
\put(21.0000,-34.0000){\makebox(0,0){$H_{6}$}}%
\put(35.0000,-34.0000){\makebox(0,0){$H_{7}$}}%
\put(49.0000,-34.0000){\makebox(0,0){$H_{8}$}}%
%
\special{sh 1.000}%
\special{ia 3500 2106 50 50  0.0000000  6.2831853}%
\special{pn 8}%
\special{pn 8}%
\special{ar 3500 2106 50 50  0.0000000  6.2831853}%
%
\special{sh 1.000}%
\special{ia 3100 2396 50 50  0.0000000  6.2831853}%
\special{pn 8}%
\special{pn 8}%
\special{ar 3100 2396 50 50  0.0000000  6.2831853}%
%
\special{sh 1.000}%
\special{ia 3900 2396 50 50  0.0000000  6.2831853}%
\special{pn 8}%
\special{pn 8}%
\special{ar 3900 2396 50 50  0.0000000  6.2831853}%
%
\special{sh 1.000}%
\special{ia 3900 2796 50 50  0.0000000  6.2831853}%
\special{pn 8}%
\special{pn 8}%
\special{ar 3900 2796 50 50  0.0000000  6.2831853}%
%
\special{sh 1.000}%
\special{ia 3100 2796 50 50  0.0000000  6.2831853}%
\special{pn 8}%
\special{pn 8}%
\special{ar 3100 2796 50 50  0.0000000  6.2831853}%
%
\special{sh 1.000}%
\special{ia 3500 3096 50 50  0.0000000  6.2831853}%
\special{pn 8}%
\special{pn 8}%
\special{ar 3500 3096 50 50  0.0000000  6.2831853}%
%
\special{pn 8}%
\special{pa 3500 3096}%
\special{pa 3100 2796}%
\special{fp}%
\special{pa 3100 2796}%
\special{pa 3100 2396}%
\special{fp}%
\special{pa 3100 2396}%
\special{pa 3500 2096}%
\special{fp}%
\special{pa 3500 2096}%
\special{pa 3900 2396}%
\special{fp}%
\special{pa 3900 2396}%
\special{pa 3900 2796}%
\special{fp}%
\special{pa 3900 2796}%
\special{pa 3500 3096}%
\special{fp}%
%
\special{sh 1.000}%
\special{ia 3200 2596 50 50  0.0000000  6.2831853}%
\special{pn 8}%
\special{pn 8}%
\special{ar 3200 2596 50 50  0.0000000  6.2831853}%
%
\special{sh 1.000}%
\special{ia 3800 2596 50 50  0.0000000  6.2831853}%
\special{pn 8}%
\special{pn 8}%
\special{ar 3800 2596 50 50  0.0000000  6.2831853}%
%
\special{pn 8}%
\special{pa 3800 2596}%
\special{pa 3900 2396}%
\special{fp}%
\special{pa 3900 2796}%
\special{pa 3800 2596}%
\special{fp}%
\special{pa 3800 2596}%
\special{pa 3200 2596}%
\special{fp}%
\special{pa 3200 2596}%
\special{pa 3100 2396}%
\special{fp}%
\special{pa 3100 2796}%
\special{pa 3200 2596}%
\special{fp}%
\special{pa 3500 3096}%
\special{pa 3500 2096}%
\special{fp}%
\put(35.0000,-19.6500){\makebox(0,0){$y_{1}$}}%
\put(36.6000,-31.0500){\makebox(0,0){$y_{2}$}}%
\put(29.3000,-23.9500){\makebox(0,0){$y_{3}$}}%
\put(29.3000,-27.9500){\makebox(0,0){$y_{4}$}}%
\put(40.5000,-23.9500){\makebox(0,0){$y_{5}$}}%
\put(40.4000,-27.9000){\makebox(0,0){$y_{6}$}}%
\put(32.9000,-24.8500){\makebox(0,0){$y_{7}$}}%
\put(37.1000,-24.8500){\makebox(0,0){$y_{8}$}}%
%
\special{sh 1.000}%
\special{ia 4900 2106 50 50  0.0000000  6.2831853}%
\special{pn 8}%
\special{pn 8}%
\special{ar 4900 2106 50 50  0.0000000  6.2831853}%
%
\special{sh 1.000}%
\special{ia 4500 2396 50 50  0.0000000  6.2831853}%
\special{pn 8}%
\special{pn 8}%
\special{ar 4500 2396 50 50  0.0000000  6.2831853}%
%
\special{sh 1.000}%
\special{ia 5300 2396 50 50  0.0000000  6.2831853}%
\special{pn 8}%
\special{pn 8}%
\special{ar 5300 2396 50 50  0.0000000  6.2831853}%
%
\special{sh 1.000}%
\special{ia 5300 2796 50 50  0.0000000  6.2831853}%
\special{pn 8}%
\special{pn 8}%
\special{ar 5300 2796 50 50  0.0000000  6.2831853}%
%
\special{sh 1.000}%
\special{ia 4500 2796 50 50  0.0000000  6.2831853}%
\special{pn 8}%
\special{pn 8}%
\special{ar 4500 2796 50 50  0.0000000  6.2831853}%
%
\special{sh 1.000}%
\special{ia 4900 3096 50 50  0.0000000  6.2831853}%
\special{pn 8}%
\special{pn 8}%
\special{ar 4900 3096 50 50  0.0000000  6.2831853}%
%
\special{pn 8}%
\special{pa 4900 3096}%
\special{pa 4500 2796}%
\special{fp}%
\special{pa 4500 2796}%
\special{pa 4500 2396}%
\special{fp}%
\special{pa 4500 2396}%
\special{pa 4900 2096}%
\special{fp}%
\special{pa 4900 2096}%
\special{pa 5300 2396}%
\special{fp}%
\special{pa 5300 2396}%
\special{pa 5300 2796}%
\special{fp}%
\special{pa 5300 2796}%
\special{pa 4900 3096}%
\special{fp}%
%
\special{sh 1.000}%
\special{ia 4600 2596 50 50  0.0000000  6.2831853}%
\special{pn 8}%
\special{pn 8}%
\special{ar 4600 2596 50 50  0.0000000  6.2831853}%
%
\special{sh 1.000}%
\special{ia 5200 2596 50 50  0.0000000  6.2831853}%
\special{pn 8}%
\special{pn 8}%
\special{ar 5200 2596 50 50  0.0000000  6.2831853}%
%
\special{pn 8}%
\special{pa 5200 2596}%
\special{pa 5300 2396}%
\special{fp}%
\special{pa 5300 2796}%
\special{pa 5200 2596}%
\special{fp}%
\special{pa 5200 2596}%
\special{pa 4600 2596}%
\special{fp}%
\special{pa 4600 2596}%
\special{pa 4500 2396}%
\special{fp}%
\special{pa 4500 2796}%
\special{pa 4600 2596}%
\special{fp}%
\put(49.0000,-19.6000){\makebox(0,0){$z_{1}$}}%
\put(50.5000,-31.0500){\makebox(0,0){$z_{2}$}}%
\put(43.3000,-23.9500){\makebox(0,0){$z_{3}$}}%
\put(43.3000,-27.9500){\makebox(0,0){$z_{4}$}}%
\put(54.4000,-23.9500){\makebox(0,0){$z_{5}$}}%
\put(54.4000,-27.9500){\makebox(0,0){$z_{6}$}}%
\put(46.9000,-24.8500){\makebox(0,0){$z_{7}$}}%
\put(51.1000,-24.8500){\makebox(0,0){$z_{8}$}}%
%
\special{sh 1.000}%
\special{ia 4900 2700 50 50  0.0000000  6.2831853}%
\special{pn 8}%
\special{pn 8}%
\special{ar 4900 2700 50 50  0.0000000  6.2831853}%
%
\special{pn 8}%
\special{pa 4900 2700}%
\special{pa 4900 3100}%
\special{fp}%
\special{pa 4900 2700}%
\special{pa 4900 2100}%
\special{fp}%
%
\special{pn 8}%
\special{pa 4900 2700}%
\special{pa 4600 2590}%
\special{fp}%
\special{pa 5200 2590}%
\special{pa 4900 2700}%
\special{fp}%
%
\special{pn 8}%
\special{pa 4900 3090}%
\special{pa 4912 3060}%
\special{pa 4926 3030}%
\special{pa 4938 3000}%
\special{pa 4950 2968}%
\special{pa 4960 2938}%
\special{pa 4972 2908}%
\special{pa 4982 2878}%
\special{pa 4990 2846}%
\special{pa 5000 2816}%
\special{pa 5006 2786}%
\special{pa 5012 2754}%
\special{pa 5018 2724}%
\special{pa 5022 2692}%
\special{pa 5024 2662}%
\special{pa 5024 2600}%
\special{pa 5022 2568}%
\special{pa 5020 2538}%
\special{pa 5016 2506}%
\special{pa 5010 2474}%
\special{pa 5004 2444}%
\special{pa 4998 2412}%
\special{pa 4982 2348}%
\special{pa 4974 2318}%
\special{pa 4964 2286}%
\special{pa 4956 2254}%
\special{pa 4946 2222}%
\special{pa 4934 2192}%
\special{pa 4914 2128}%
\special{pa 4902 2096}%
\special{pa 4900 2090}%
\special{sp}%
\put(47.8000,-27.3000){\makebox(0,0){$z_{9}$}}%
%
\special{pn 8}%
\special{pa 3100 2400}%
\special{pa 3900 2400}%
\special{fp}%
\special{pa 3900 2800}%
\special{pa 3100 2800}%
\special{fp}%
%
\special{pn 8}%
\special{pa 4500 2800}%
\special{pa 5300 2800}%
\special{fp}%
\special{pa 5300 2400}%
\special{pa 4500 2400}%
\special{fp}%
\end{picture}%

\caption{Graphs $H_{i}$}
\label{f1}
\end{center}
\end{figure}

We first give a characterization of graphs as in (F1).
A generalized comb is {\it pointed} if all of its leaf-cliques consist of exactly one vertex.
Let $\H_{0}$ be the family of pointed generalized combs.
For each $i~(1\leq i\leq 8)$, let $H_{i}$ be the graph depicted in Figure~\ref{f1}.
For each integer $i~(1\leq i\leq 5)$, 
the vertices of $H_{i}$ enclosed with a circle are called {\it expandable vertices}.
Also, for an expandable vertex $a$ of $H_{i}$, 
{\it expanding} of $a$ to a clique $C$ is the operation 
replacing $a$ to $C$ and adding additional edges between $u\in V(H_{i})-\{a\}$ and $C$ 
if $au \in E(H_{i})$. 
Let $U_{i}$ be the set of expandable vertices of $H_{i}$.
For a family $\C=\{C_{a}\mid a\in U_{i}\}$ of vertex-disjoint cliques indexed by $a$, 
the graph $H_{i}(\C)$ is obtained from $H_{i}$ by expanding each vertex $a\in U_{i}$ to the clique $C_{a}$ (see Figure~\ref{f2}).
\begin{figure}[htb]
\begin{center}
\unitlength 0.1in
\begin{picture}( 35.9000,  9.6000)(  9.1000,-16.3500)
%
\special{sh 1.000}%
\special{ia 1200 1400 50 50  0.0000000  6.2831853}%
\special{pn 8}%
\special{pn 8}%
\special{ar 1200 1400 50 50  0.0000000  6.2831853}%
%
\special{sh 1.000}%
\special{ia 1600 1000 50 50  0.0000000  6.2831853}%
\special{pn 8}%
\special{pn 8}%
\special{ar 1600 1000 50 50  0.0000000  6.2831853}%
%
\special{sh 1.000}%
\special{ia 2000 1400 50 50  0.0000000  6.2831853}%
\special{pn 8}%
\special{pn 8}%
\special{ar 2000 1400 50 50  0.0000000  6.2831853}%
%
\special{pn 8}%
\special{ar 2000 1400 100 100  0.0000000  6.2831853}%
%
\special{pn 8}%
\special{ar 1200 1400 100 100  0.0000000  6.2831853}%
%
\special{pn 8}%
\special{pa 1200 1400}%
\special{pa 1600 1000}%
\special{fp}%
\special{pa 1600 1000}%
\special{pa 2000 1400}%
\special{fp}%
\special{pa 2000 1400}%
\special{pa 1200 1400}%
\special{fp}%
%
\special{pn 8}%
\special{pa 1400 1600}%
\special{pa 1280 1480}%
\special{fp}%
\special{sh 1}%
\special{pa 1280 1480}%
\special{pa 1314 1542}%
\special{pa 1318 1518}%
\special{pa 1342 1514}%
\special{pa 1280 1480}%
\special{fp}%
\special{pa 1480 1600}%
\special{pa 1910 1470}%
\special{fp}%
\special{sh 1}%
\special{pa 1910 1470}%
\special{pa 1840 1470}%
\special{pa 1860 1486}%
\special{pa 1852 1508}%
\special{pa 1910 1470}%
\special{fp}%
\put(14.4000,-17.0000){\makebox(0,0){extendable}}%
\put(20.4000,-7.4000){\makebox(0,0){non-extendable}}%
%
\special{pn 20}%
\special{pa 2400 1200}%
\special{pa 3000 1200}%
\special{fp}%
\special{sh 1}%
\special{pa 3000 1200}%
\special{pa 2934 1180}%
\special{pa 2948 1200}%
\special{pa 2934 1220}%
\special{pa 3000 1200}%
\special{fp}%
%
\special{sh 1.000}%
\special{ia 4000 1000 50 50  0.0000000  6.2831853}%
\special{pn 8}%
\special{pn 8}%
\special{ar 4000 1000 50 50  0.0000000  6.2831853}%
%
\special{pn 8}%
\special{ar 3600 1400 180 180  0.0000000  6.2831853}%
%
\special{pn 8}%
\special{ar 4400 1400 100 100  0.0000000  6.2831853}%
%
\special{pn 8}%
\special{pa 4000 1000}%
\special{pa 3500 1250}%
\special{fp}%
\special{pa 4000 1000}%
\special{pa 3770 1480}%
\special{fp}%
%
\special{pn 8}%
\special{pa 4000 1000}%
\special{pa 4450 1310}%
\special{fp}%
\special{pa 4000 1000}%
\special{pa 4330 1480}%
\special{fp}%
%
\special{pn 8}%
\special{pa 3650 1230}%
\special{pa 4430 1300}%
\special{fp}%
%
\special{pn 8}%
\special{pa 4430 1500}%
\special{pa 3620 1580}%
\special{fp}%
\put(40.0000,-14.0000){\makebox(0,0){$+$}}%
\put(38.2000,-11.7000){\makebox(0,0){$+$}}%
\put(42.1000,-12.1000){\makebox(0,0){$+$}}%
\put(17.0000,-10.1000){\makebox(0,0){$u$}}%
\put(10.0000,-14.1000){\makebox(0,0){$v$}}%
\put(22.0000,-14.1000){\makebox(0,0){$w$}}%
\put(46.0000,-14.1000){\makebox(0,0){$C_{w}$}}%
\put(32.9000,-14.0000){\makebox(0,0){$C_{v}$}}%
\put(41.5000,-9.1000){\makebox(0,0){$u$}}%
%
\special{pn 8}%
\special{pa 1750 820}%
\special{pa 1650 940}%
\special{fp}%
\special{sh 1}%
\special{pa 1650 940}%
\special{pa 1708 902}%
\special{pa 1684 900}%
\special{pa 1678 876}%
\special{pa 1650 940}%
\special{fp}%
\end{picture}%
\caption{Expanding vertices to cliques}
\label{f2}
\end{center}
\end{figure}
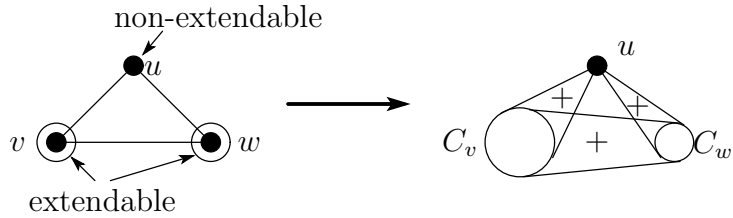
Let
$$
\H_{i}=\{H_{i}(\C)\mid \C=\{C_{a}\mid a\in U_{i}\}\mbox{ is a family of vertex-disjoint cliques indexed by }a\}.
$$
Note that $H_{i}\in \H_{i}$.
For each $j~(6\leq j\leq 8)$, let $\H_{j}=\{H_{j}\}$.

\begin{thm}
\label{thm1z}
Let $G$ be a connected $\{K_{1,3},Z_{2}\}$-free graph.
Then $G$ is not $B_{1,1}$-free if and only if $G\in \bigcup _{0\leq i\leq 8}\H_{i}$.
\end{thm}

We next consider giving a characterization of graphs as depicted in (F2) and (F3).
Let $l\geq 5$ be an integer, and let $L_{0},L_{1},\ldots ,L_{l}$ be vertex-disjoint cliques.
The graph $F_{p}=F_{p}(L_{1},\ldots ,L_{l})$ is obtained from $L_{1}\cup \cdots \cup L_{l}$ by joining every vertex of $L_{i}$ to all vertices of $L_{i+1}$ for $1\leq i\leq l-1$, and we call $F_{p}$ a {\it fat $l$-path} (or simply a {\it fat path}).
In this context, $L_{i}~(1\leq i\leq l)$ are called {\it fundamental cliques} of $F_{p}$.
The graph $F_{c}=F_{c}(L_{0},\ldots ,L_{l})$ is obtained from $L_{0}\cup \cdots \cup L_{l}$ by joining every vertex of $L_{i}$ to all vertices of $L_{i+1}$ for $0\leq i\leq l$ where the indices are calculated modulo $l+1$, and we call $F_{c}$ a {\it fat $l$-cycle} (or simply a {\it fat cycle}).
In this context, $L_{i}~(0\leq i\leq l)$ are called {\it fundamental cliques} of $F_{c}$.
Note that fat $l$-paths have $l$ fundamental cliques but fat $l$-cycles have $l+1$ fundamental cliques.
Let $\P(l)$ be the family of fat $i$-paths and fat $i$-cycles for all $i\geq l$.

We give the following characterization.

\begin{thm}
\label{thm1}
Let $G$ be a connected $\{K_{1,3},B_{1,1}\}$-free graph.
Then $G$ is not $P_{5}$-free if and only if $G\in \P(5)$.
\end{thm}

\begin{thm}
\label{thm2}
Let $G$ be a connected $\{K_{1,3},B_{1,2}\}$-free graph.
Then $G$ is not $P_{6}$-free if and only if $G\in \P(6)$.
\end{thm}

We prove the following generalization of Theorems~\ref{thm1} and \ref{thm2}.

\begin{thm}
\label{thm3}
Let $m\geq 1$ be an integer, and let $G$ be a connected $\{K_{1,3},B_{1,m}\}$-free graph.
Then $G$ is not $P_{\max\{3m,m+4\}}$-free if and only if $G\in \P(\max\{3m,m+4\})$.
\end{thm}

\begin{remark}
There are infinitely many connected $\{K_{1,3},B_{1,m}\}$-free but not $P_{\max\{3m-1,m+3\}}$-free graphs which are neither fat paths nor fat cycles:
Fix an integer $m\geq 1$.
Let $F=F_{p}(L_{1},\ldots ,L_{\max\{3m-1,m+3\}})$ be a fat path, and let $K$ be a clique with $V(F)\cap K=\emptyset $.
Let $F'$ be the graph obtained from $F\cup K$ by joining each vertex of $K$ to each vertex of $L_{\max\{m-1,1\}}\cup L_{\max\{m,2\}}\cup L_{\max\{2m,3\}}\cup L_{\max\{2m+1,4\}}$ (see Figure~\ref{notpath-free}).
Then we see that $F'$ is a connected $\{K_{1,3},B_{1,m}\}$-free but not $P_{\max\{3m-1,m+3\}}$-free graph.
Therefore the order of the path in Theorem~\ref{thm3} is best possible if we require the targets to be graphs with a simple structure.
\end{remark}

\begin{figure}
\begin{center}
{\unitlength 0.1in
\begin{picture}( 55.0000, 10.1500)(  4.0500,-17.4000)
%
\special{pn 8}%
\special{ar 600 1400 100 100  0.0000000  6.2831853}%
%
\special{pn 8}%
\special{ar 900 1400 100 100  0.0000000  6.2831853}%
%
\special{pn 8}%
\special{ar 1500 1400 100 100  0.0000000  6.2831853}%
%
\special{pn 8}%
\special{ar 1200 1400 100 100  0.0000000  6.2831853}%
%
\special{pn 20}%
\special{pa 700 1400}%
\special{pa 800 1400}%
\special{fp}%
%
\special{pn 20}%
\special{pa 1000 1400}%
\special{pa 1100 1400}%
\special{fp}%
\special{pa 1300 1400}%
\special{pa 1400 1400}%
\special{fp}%
\put(6.0000,-16.0000){\makebox(0,0){$L_{1}$}}%
%
\special{pn 8}%
\special{ar 1050 1000 100 100  0.0000000  6.2831853}%
%
\special{pn 20}%
\special{pa 600 1300}%
\special{pa 970 1060}%
\special{fp}%
%
\special{pn 20}%
\special{pa 1130 1060}%
\special{pa 1500 1300}%
\special{fp}%
%
\special{pn 20}%
\special{pa 900 1300}%
\special{pa 1000 1090}%
\special{fp}%
\special{pa 1100 1090}%
\special{pa 1200 1300}%
\special{fp}%
\put(10.5000,-8.0000){\makebox(0,0){$K$}}%
\put(10.5000,-18.0000){\makebox(0,0){$(m=1)$}}%
\put(12.0000,-16.0000){\makebox(0,0){$L_{3}$}}%
\put(9.0000,-16.0000){\makebox(0,0){$L_{2}$}}%
\put(15.0000,-16.0000){\makebox(0,0){$L_{4}$}}%
%
\special{pn 8}%
\special{ar 2196 1396 100 100  0.0000000  6.2831853}%
%
\special{pn 8}%
\special{ar 2796 1396 100 100  0.0000000  6.2831853}%
%
\special{pn 8}%
\special{ar 3096 1396 100 100  0.0000000  6.2831853}%
%
\special{pn 8}%
\special{ar 3696 1396 100 100  0.0000000  6.2831853}%
%
\special{pn 8}%
\special{ar 5806 1396 100 100  1.6704650  1.5707963}%
%
\special{pn 4}%
\special{sh 1}%
\special{ar 2496 1396 8 8 0  6.28318530717959E+0000}%
\special{sh 1}%
\special{ar 2406 1396 8 8 0  6.28318530717959E+0000}%
\special{sh 1}%
\special{ar 2606 1396 8 8 0  6.28318530717959E+0000}%
\special{sh 1}%
\special{ar 2606 1396 8 8 0  6.28318530717959E+0000}%
%
\special{pn 4}%
\special{sh 1}%
\special{ar 3996 1396 8 8 0  6.28318530717959E+0000}%
\special{sh 1}%
\special{ar 3906 1396 8 8 0  6.28318530717959E+0000}%
\special{sh 1}%
\special{ar 4106 1396 8 8 0  6.28318530717959E+0000}%
\special{sh 1}%
\special{ar 4106 1396 8 8 0  6.28318530717959E+0000}%
%
\special{pn 4}%
\special{sh 1}%
\special{ar 5496 1396 8 8 0  6.28318530717959E+0000}%
\special{sh 1}%
\special{ar 5406 1396 8 8 0  6.28318530717959E+0000}%
\special{sh 1}%
\special{ar 5606 1396 8 8 0  6.28318530717959E+0000}%
\special{sh 1}%
\special{ar 5606 1396 8 8 0  6.28318530717959E+0000}%
%
\special{pn 8}%
\special{ar 4000 1000 100 100  0.0000000  6.2831853}%
%
\special{pn 8}%
\special{ar 3396 1396 100 100  0.0000000  6.2831853}%
%
\special{pn 20}%
\special{pa 2896 1396}%
\special{pa 2996 1396}%
\special{fp}%
%
\special{pn 20}%
\special{pa 3196 1396}%
\special{pa 3296 1396}%
\special{fp}%
\special{pa 3496 1396}%
\special{pa 3596 1396}%
\special{fp}%
%
\special{pn 20}%
\special{pa 2696 1396}%
\special{pa 2646 1396}%
\special{fp}%
%
\special{pn 20}%
\special{pa 2346 1396}%
\special{pa 2296 1396}%
\special{fp}%
%
\special{pn 20}%
\special{pa 3796 1396}%
\special{pa 3846 1396}%
\special{fp}%
%
\special{pn 8}%
\special{ar 4296 1396 100 100  0.0000000  6.2831853}%
%
\special{pn 8}%
\special{ar 4596 1396 100 100  0.0000000  6.2831853}%
%
\special{pn 8}%
\special{ar 5196 1396 100 100  0.0000000  6.2831853}%
%
\special{pn 8}%
\special{ar 4896 1396 100 100  0.0000000  6.2831853}%
%
\special{pn 20}%
\special{pa 4396 1396}%
\special{pa 4496 1396}%
\special{fp}%
%
\special{pn 20}%
\special{pa 4696 1396}%
\special{pa 4796 1396}%
\special{fp}%
\special{pa 4996 1396}%
\special{pa 5096 1396}%
\special{fp}%
%
\special{pn 20}%
\special{pa 4196 1396}%
\special{pa 4146 1396}%
\special{fp}%
%
\special{pn 20}%
\special{pa 5296 1396}%
\special{pa 5346 1396}%
\special{fp}%
%
\special{pn 20}%
\special{pa 5696 1396}%
\special{pa 5646 1396}%
\special{fp}%
\put(39.9000,-7.9000){\makebox(0,0){$K$}}%
\put(21.9500,-15.9500){\makebox(0,0){$L_{1}$}}%
\put(30.4500,-15.9500){\makebox(0,0){$L_{m-1}$}}%
\put(34.4500,-15.9500){\makebox(0,0){$L_{m}$}}%
\put(45.4500,-15.9500){\makebox(0,0){$L_{2m}$}}%
\put(49.4500,-15.9500){\makebox(0,0){$L_{2m+1}$}}%
\put(57.9500,-15.9500){\makebox(0,0){$L_{3m-1}$}}%
\put(40.0500,-18.0500){\makebox(0,0){$(m\geq 2)$}}%
%
\special{pn 20}%
\special{pa 3096 1306}%
\special{pa 3900 1040}%
\special{fp}%
\special{pa 4090 1040}%
\special{pa 4900 1300}%
\special{fp}%
%
\special{pn 20}%
\special{pa 3400 1300}%
\special{pa 3940 1070}%
\special{fp}%
\special{pa 4060 1070}%
\special{pa 4600 1300}%
\special{fp}%
\end{picture}}%
\caption{Graph $F'$}
\label{notpath-free}
\end{center}
\end{figure}

We prove Theorems~\ref{thm1z} and \ref{thm3} in Subsections~\ref{subsec2.1} and \ref{subsec2.2}, respectively.

\subsection{Applications}\label{subsec1.3}

Duffus, Jacobson and Gould~\cite{DJG} proved that every $2$-connected $\{K_{1,3},N\}$-free graph has a Hamiltonian cycle.
Since we can verify that every $2$-connected generalized comb has a Hamiltonian cycle, this together with Theorem~\ref{ThmA} implies that every $2$-connected $\{K_{1,3},B_{1,2}\}$-free graph has a Hamiltonian cycle.
In other words, Duffus-Jacobson-Gould theorem and Theorem~\ref{ThmA} provide an alternative proof of Theorem~\ref{Be} in Subsection~\ref{subsec1.3.1}.
Our main results have similar applications.

For our argument, we introduce a notation related to forbidden subgraphs.
For two families $\H_{1}$ and $\H_{2}$ of forbidden subgraphs, we write $\H_{1}\leq \H_{2}$ if for every $H_{2}\in \H_{2}$, there exists $H_{1}\in \H_{1}$ such that $H_{1}$ is an induced subgraph of $H_{2}$.
Note that if $\H_{1}\leq \H_{2}$, then every $\H_{1}$-free graph is also $\H_{2}$-free.

\subsubsection{Hamiltonian cycles}\label{subsec1.3.1}

In the study of forbidden subgraphs, it is a fundamental problem to characterize the forbidden pairs assuring some properties $P$.
When we consider such problems, we often assume a trivial necessary condition of $P$ (for example, when we consider the existence of a Hamiltonian cycle, it is natural to assume the $2$-connectedness).
Bedrossian~\cite{Be} characterized the forbidden pairs for the existence of a Hamiltonian cycle as follows:

\begin{Thm}[Bedrossian~\cite{Be}]
\label{Be2}
Let $\H$ be a forbidden pair. 
Then every $2$-connected $\H$-free graph has a Hamiltonian cycle if and only if either $\H \le \{K_{1,3},N\}$ or $\H \le \{K_{1,3},B_{1,2}\}$ or $\H \le \{K_{1,3},P_{6}\}$.
\end{Thm}

Bedrossian's characterization depends on, for example, the following theorems.

\begin{Thm}[Broersma and Veldman~\cite{BV}]
\label{BV}
Every $2$-connected $\{K_{1,3},P_{6}\}$-free graph has a Hamiltonian cycle.
\end{Thm}

\begin{Thm}[Bedrossian~\cite{Be}]
\label{Be}
Every $2$-connected $\{K_{1,3},B_{1,2}\}$-free graph has a Hamiltonian cycle. 
\end{Thm}

Since any $2$-connected fat $i$-paths and any $2$-connected fat $i$-cycles have a Hamiltonian cycle for $i\geq 6$, Theorems~\ref{thm2} and \ref{BV} give an alternative proof of Theorem~\ref{Be}.

\subsubsection{Halin graphs}\label{subsec1.3.2}

A graph is {\it planar} if it can be embedded in the plane without edge-crossing, and such an embedded graph is called a {\it plane graph}.
A {\it Halin graph}, defined by Halin~\cite{H}, is a plane graph consisting of a tree $T$ without vertices of degree $2$ and a cycle $C$ induced by the leaves of $T$ (and we often write a Halin graph $H$ as $H=T\cup C$).
If a graph $G$ contains a Halin graph as a spanning subgraph, then it is called a {\it spanning Halin subgraph} of $G$.
In~\cite{CHOST}, the following conjecture was proposed.

\begin{con}[Chen, Han, O, Shan and Tsuchiya~\cite{CHOST}]
\label{conj1}
Let $\H$ be a forbidden pair. 
Then every $3$-connected $\H$-free graph has a spanning Halin subgraph if and only if either $\H \le \{K_{1,3},Z_{3}\}$ or $\H \le \{K_{1,3},B_{1,2}\}$.
\end{con}

The ``only if'' part of Conjecture~\ref{conj1} was already proved in \cite{CHOST}.
Also, as a partial answer for ``if'' part of the conjecture, the following theorem was proved. 

\begin{Thm}[Chen, Han, O, Shan and Tsuchiya~\cite{CHOST}]
\label{thm:P5}
Every $3$-connected $\{K_{1,3},P_{5}\}$-free graph has a spanning Halin subgraph. 
\end{Thm}

As corollaries of Theorems~\ref{thm1z}, \ref{thm1} and~\ref{thm:P5}, we obtain other partial answers for ``if'' part of Conjecture~\ref{conj1} (and we give those detail in Appendix).

\begin{thm}
\label{thm:SHS2}
Every $3$-connected $\{K_{1,3},B_{1,1}\}$-free graph has a spanning Halin subgraph. 
\end{thm}

\begin{thm}
\label{thm:SHS3}
Every $3$-connected $\{K_{1,3},Z_{2}\}$-free graph has a spanning Halin subgraph. 
\end{thm}

\subsubsection{Independence numbers}\label{subsec1.3.3}

The {\it independence number} of a graph $G$ is the maximum cardinality of an independent set of $G$.
In~\cite{BH}, Brandst\"{a}dt and Hammer found a polynomial-time algorithm for determining the independence number of $\{K_{1,3},P_{5}\}$-free graphs.

Let $H$ be a graph belonging to $\P(5)$, and let $Q$ be an induced path of $H$ having three vertices.
Then any maximal induced paths and any maximal induced cycles containing $Q$ pass through each fundamental cliques of $H$ exactly once.
By using the fact above, for a given graph $G$, we can decide whether $G$ belongs to $\P(5)$ or not in polynomial-time (and we omit its precise algorithm).
This together with Theorem~\ref{thm1} assures that we can determine the independence number of $\{K_{1,3},B_{1,1}\}$-free graphs in polynomial-time.

\section{Proof of main results}\label{sec2}

In this section, we prove Theorems~\ref{thm1z} and \ref{thm3}.

\subsection{Proof of Theorem~\ref{thm1z}}\label{subsec2.1}

\begin{lem}
\label{lem2.1}
Let $G$ be a connected $\{K_{1,3},Z_{2}\}$-free graph which contains an induced subgraph $N$.
Then $G$ is a pointed generalized comb.
\end{lem}
\proof
Since $G$ is $Z_{2}$-free and $Z_2$ is an induced subgraph of $B_{1,2}$, $G$ is also $B_{1,2}$-free.
This, together with Theorem~\ref{ThmA}, implies that $G$ is a generalized comb.
We only show that every leaf-clique consists of a single vertex.
Let $L_{i}~(1\leq i\leq m)$ be the leaf-cliques of $G$, and let $R_{i}$ be the root of $L_{i}$.
On the contrary, we may assume that $|L_{1}|\geq 2$.
Let $a_{1},a_{2}\in L_{1}$ with $a_{1}\not=a_{2}$, $a_{3}\in R_{1}$, $a_{4}\in R_{2}$ and $a_{5}\in L_{2}$.
Then $G[\{a_{1},a_{2},a_{3},a_{4},a_{5}\}]\cong Z_{2}$, giving a contradiction.
Hence $G$ is a pointed generalized comb.
\qed

In the following lemmas (Lemmas~\ref{lem2.2}--\ref{lem2.8}), we follow the labels given in Figures~\ref{f1} and~\ref{f2}.

\begin{lem}
\label{lem2.2}
Let $G$ be a connected $\{K_{1,3},Z_{2},N\}$-free graph which contains an induced subgraph $H=H_{1}(\{C_{s_{3}},C_{s_{4}},C_{s_{5}}\})$, where $C_{s_{3}}$, $C_{s_{4}}$ and $C_{s_{5}}$ are vertex-disjoint cliques.
Then for each vertex $a\in V(G)-V(H)$ with $N_{G}(a)\cap V(H)\not=\emptyset $, one of the following holds:
\begin{enumerate}[{\upshape(i)}]
\item
$G[V(H)\cup \{a\}]\in \H_{1}$,
\item
for some $i\in \{1,2\}$, $N_{G}(a)\cap V(H)=\{s_{i}\}\cup C_{s_{5}}$ and $|C_{s_{5-i}}|=1$ $($and so $G[V(H)\cup \{a\}]\in \H_{2})$ ,
\item
for some $i\in \{3,4\}$, $N_{G}(a)\cap V(H)=\{s_{1},s_{2}\}\cup C_{s_{i}}$ and $|C_{s_{7-i}}|=|C_{s_{5}}|=1$ $($and so $G[V(H)\cup \{a\}]\in \H_{3})$, or
\item
$N_{G}(a)\cap V(H)=\{s_{1},s_{2}\}\cup C_{s_{3}}\cup C_{s_{4}}$ and $|C_{s_{5}}|=1$ $($and so $G[V(H)\cup \{a\}]\in \H_{4})$.
\end{enumerate}
\end{lem}
\proof
For each $i\in \{3,4,5\}$, we take a vertex $b_{i}$ as follows:
If $N_{G}(a)\cap C_{s_{i}} \neq \emptyset$, let $b_{i}\in N_{G}(a)\cap C_{s_{i}}$; otherwise (i.e., $N_{G}(a)\cap C_{s_{i}}=\emptyset$), let $b_{i} \in C_{s_{i}}$.

\medskip
\noindent
\textbf{Case 1:} $N_{G}(a) \cap C_{s_{5}} \neq \emptyset$.

If $as_{1},as_{2}\in E(G)$, then $G[\{a,s_{1},s_{2},b_{5}\}]\cong K_{1,3}$, giving a contradiction.
Thus $as_{1}\not\in E(G)$ or $as_{2}\not\in E(G)$.
We may assume that $as_{1}\not\in E(G)$.

If $N_G(a) \cap C_{s_{3}} \neq \emptyset$ and $ab\not\in E(G)$ for some $b\in C_{s_{4}}$, then $G[\{b_{3},a,b,s_{1}\}]\cong K_{1,3}$; if $N_G(a) \cap C_{s_{4}} \neq \emptyset$ and $ab\not\in E(G)$ for some $b\in C_{s_{3}}$, then $as_{2}\in E(G)$ because $G[\{b_{4},a,b,s_{2}\}]\not\cong K_{1,3}$, and hence $G[\{a,s_{2},b_{4},b,s_{1}\}]\cong Z_{2}$.
In either case, we get a contradiction.
This implies that either $C_{s_{3}}\cup C_{s_{4}}\subseteq N_{G}(a)$ or $N_{G}(a)\cap (C_{s_{3}}\cup C_{s_{4}})=\emptyset $.

\medskip
\noindent
\textbf{Subcase 1.1:} $C_{s_{3}}\cup C_{s_{4}}\subseteq N_{G}(a)$.

If $ab\not\in E(G)$ for some $b\in C_{s_{5}}$, then $G[\{b_{3},a,s_{1},b\}]\cong K_{1,3}$, giving a contradiction.
Thus $C_{s_{5}}\subseteq N_{G}(a)$.
If $as_{2}\in E(G)$, let $C'_{s_{i}}=C_{s_{i}}~(i\in \{3,5\})$ and $C'_{s_{4}}=C_{s_{4}}\cup \{a\}$; if $as_{2}\not \in E(G)$, let $C'_{s_{i}}=C_{s_{i}}~(i\in \{3,4\})$ and $C'_{s_{5}}=C_{s_{5}}\cup \{a\}$.
Then $G[V(H)\cup \{a\}]=H_{1}(\{C'_{s_{3}},C'_{s_{4}},C'_{s_{5}}\})\in \H_{1}$.

\medskip
\noindent
\textbf{Subcase 1.2:} $N_{G}(a)\cap (C_{s_{3}}\cup C_{s_{4}})=\emptyset $.

Since $G[\{a,s_{1},s_{2},b_{5},b_{3},b_{4}\}]\not\cong N$, we have $as_{2}\in E(G)$.
If $ab\not\in E(G)$ for some $b\in C_{s_{5}}$, then $G[\{b,b_{3},b_{4},s_{2},a\}]\cong Z_{2}$, giving a contradiction.
Thus $C_{s_{5}}\subseteq N_{G}(a)$, and hence $N_{G}(a)\cap V(H)=\{s_{2}\}\cup C_{s_{5}}$.
If $|C_{s_{3}}|\geq 2$, then $G[\{b_{3},b,b_{4},s_{2},a\}]\cong Z_{2}$ where $b\in C_{s_{3}}-\{b_{3}\}$, giving a contradiction.
Thus $|C_{s_{3}}|=1$.

\medskip
\noindent
\textbf{Case 2:} $N_{G}(a)\cap C_{s_{5}}=\emptyset $ (i.e., $ab_{5}\not\in E(G)$).

\begin{claim}
\label{cl2.2.1}
For each $i\in \{3,4\}$, if $N_{G}(a) \cap C_{s_{i}}\neq\emptyset$, then $N_{G}(a)\supseteq \{s_{1},s_{2}\}\cup C_{s_{i}}$.
\end{claim}
\proof
We may assume $i=3$.
Since $G[\{b_{3},a,b_{5},s_{1}\}]\not\cong K_{1,3}$, we have $as_{1}\in E(G)$.
By the same argument, if $N_{G}(a) \cap C_{s_{4}}\neq\emptyset$,
then $as_{2}\in E(G)$.
Since $G[\{a,s_{1},b_{3},b_{4},s_{2}\}]\not\cong Z_{2}$, we have $as_{2}\in E(G)$ or $ab_{4}\in E(G)$.
In either case, we have $as_{2}\in E(G)$.
If $ab\not\in E(G)$ for some $b\in C_{s_{3}}$, then $G[\{b_{5},b,b_{3},a,s_{2}\}]\cong Z_{2}$, giving a contradiction.
Thus $C_{s_{3}}\subseteq N_{G}(a)$.
\qed

Suppose $N_{G}(a)\cap (C_{s_{3}}\cup C_{s_{4}})=\emptyset $.
Since $N_{G}(a)\cap V(H)\not=\emptyset $, we have $as_{i}\in E(G)$ for some $i\in \{1,2\}$.
Hence $G[\{b_{5},b_{5-i},b_{i+2},s_{i},a\}]\cong Z_{2}$, giving a contradiction.
Thus $N_{G}(a)\cap (C_{s_{3}}\cup C_{s_{4}})\not=\emptyset $.
We may assume that $N_{G}(a) \cap C_{s_{3}}\neq\emptyset$.
This together with Claim~\ref{cl2.2.1} forces $\{s_{1},s_{2}\}\cup C_{s_{3}}\subseteq N_{G}(a)$.
If $|C_{s_{5}}|\geq 2$, then $G[\{b_{5},b,b_{3},a,s_{2}\}]\cong Z_{2}$ where $b\in C_{s_{5}}-\{b_{5}\}$, giving a contradiction.
Thus $|C_{s_{5}}|=1$.

If $N_G(a) \cap C_{s_{4}}\neq\emptyset$, then $C_{s_{4}}\subseteq N_{G}(a)$ by Claim~\ref{cl2.2.1}, and hence (iv) holds.
Thus we may assume that $N_G(a) \cap C_{s_{i}}=\emptyset$ (i.e., $N_{G}(a)\cap V(H)=\{s_{1},s_{2}\}\cup C_{s_{3}}$).
If $|C_{s_{4}}|\geq 2$, then $G[\{b_{4},b,s_{2},a,s_{1}\}]\cong Z_{2}$ in $G$ where $b\in C_{s_{4}}-\{b_{4}\}$, giving a contradiction.
Hence $|C_{s_{4}}|=1$, and so (iii) holds.
\qed

\begin{lem}
\label{lem2.3}
Let $G$ be a connected $\{K_{1,3},Z_{2},N\}$-free graph which contains an induced subgraph $H=H_{2}(\{C_{t_{3}},C_{t_{4}}\})$, where $C_{t_{3}}$ and $C_{t_{4}}$ are vertex-disjoint cliques.
Then for each vertex $a\in V(G)-V(H)$ with $N_{G}(a)\cap V(H)\not=\emptyset $, one of the following holds:
\begin{enumerate}[{\upshape(i)}]
\item
$G[V(H)\cup \{a\}]\in \H_{2}$,
\item
for some $i\in \{3,4\}$, $N_{G}(a)\cap V(H)=\{t_{1},t_{i+2}\}\cup C_{t_{i}}$ and $|C_{t_{7-i}}|=1$ (and so $G[V(H)\cup \{a\}]\in \H_{5}$), or
\item
$N_{G}(a)\cap V(H)=\{t_{1},t_{2},t_{5},t_{6}\}$ and $|C_{t_{3}}|=|C_{t_{4}}|=1$ (and so $G[V(H)\cup \{a\}]\in \H_{6}$).
\end{enumerate}
\end{lem}
\proof
For each $i\in \{3,4\}$, let $b_{i}\in C_{t_{i}}$.
For each $i\in \{5,6\}$, we note that the graph $B_{i}:=H-t_{i}$ belongs to $\H_{1}$.

\medskip
\noindent
\textbf{Case 1:} $N_{G}(a)\cap (C_{t_{3}}\cup C_{t_{4}})=\emptyset $.

Since $N_{G}(a)\cap V(H)\not=\emptyset $, $N_{G}(a)\cap V(B_{i})\not=\emptyset $ for some $i\in \{5,6\}$.
We may assume that $N_{G}(a)\cap V(B_{5})\not=\emptyset $.
Since $N_{G}(a)\cap (C_{t_{3}}\cup C_{t_{4}})=\emptyset $, we have $N_{G}(a)\cap V(B_{5})=\{t_{1},t_{2},t_{6}\}$ and $|C_{t_{4}}|=1$ by Lemma~\ref{lem2.2}.
In particular, $N_{G}(a)\cap V(B_{6})\not=\emptyset $.
Then again by Lemma~\ref{lem2.2}, $N_{G}(a)\cap V(B_{6})=\{t_{1},t_{2},t_{5}\}$ and $|C_{t_{3}}|=1$.
This implies that $N_{G}(a)\cap V(H)=\{t_{1},t_{2},t_{5},t_{6}\}$ and $|C_{t_{3}}|=|C_{t_{4}}|=1$, as desired.

\medskip
\noindent
\textbf{Case 2:} $N_{G}(a)\cap (C_{t_{3}}\cup C_{t_{4}})\not=\emptyset $.

We may assume that $N_{G}(a)\cap C_{t_{3}}\not=\emptyset $.
If $N_{G}(a)\cap V(B_{5})=\{t_{6}\}\cup C_{t_{3}}$, then either $N_{G}(a)\cap V(B_{6})=C_{t_{3}}$ or $N_{G}(a)\cap V(B_{6})=\{t_{5}\}\cup C_{t_{3}}$, which contradicts Lemma~\ref{lem2.2}.
Thus, by Lemma~\ref{lem2.2}, we have either $G[V(B_{5})\cup \{a\}]\in \H_{1}$, or $N_{G}(a)\cap V(B_{5})=\{t_{1}\}\cup C_{t_{3}}$ and $|C_{t_{4}}|=1$.

\medskip
\noindent
\textbf{Subcase 2.1:} $G[V(B_{5})\cup \{a\}]\in \H_{1}$.

We see that $\{t_{2}\}\cup C_{t_{3}}\cup C_{t_{4}}\subseteq N_{G}(a)$.
Since $G[\{a,t_{2},b_{4},t_{6},t_{5}\}]\not\cong Z_{2}$, we have $at_{5}\in E(G)$ or $at_{6}\in E(G)$.
We may assume that $at_{5}\in E(G)$.
If $at_{1}\in E(G)$, then $G[\{a,t_{1},b_{4},t_{5}\}]\cong K_{1,3}$, giving a contradiction.
Thus $at_{1}\not\in E(G)$. 
So, $at_{6}\not\in E(G)$ because $G[\{t_{5},t_{6},a,t_{2},t_{1}\}]\not\cong Z_{2}$.
Hence we get $N_{G}(a)\cap V(H)=\{t_{2},t_{5}\}\cup C_{t_{3}}\cup C_{t_{4}}$.
Consequently, $G[V(H)\cup \{a\}]=H_{2}(\{C'_{t_{3}},C'_{t_{4}}\})$ where $C'_{s_{3}}=C_{s_{3}}\cup \{a\}$ and $C'_{s_{4}}=C_{s_{4}}$, as desired.

\medskip
\noindent
\textbf{Subcase 2.2:} $N_{G}(a)\cap V(B_{5})=\{t_{1}\}\cup C_{t_{3}}$ and $|C_{t_{4}}|=1$.

Since $G[\{b_{3},a,t_{2},t_{5}\}]\not\cong K_{1,3}$, we have $at_{5}\in E(G)$.
Hence $N_{G}(a)\cap V(H)=\{t_{1},t_{5}\}\cup C_{t_{3}}$ and $|C_{t_{4}}|=1$.
\qed

\begin{lem}
\label{lem2.4}
Let $G$ be a connected $\{K_{1,3},Z_{2},N\}$-free graph which contains an induced subgraph $H=H_{3}(\{C_{u_{6}}\})$, where $C_{u_{6}}$ is a clique.
Then for each vertex $a\in V(G)-V(H)$ with $N_{G}(a)\cap V(H)\not=\emptyset $, one of the following holds:
\begin{enumerate}[{\upshape(i)}]
\item
$G[V(H)\cup \{a\}]\in \H_{3}$,
\item
$N_{G}(a)\cap V(H)=\{u_{1},u_{i},u_{7-i}\}$ for some $i\in \{2,3\}$ and $|C_{u_{6}}|=1$ (and so $G[V(H)\cup \{a\}]\in \H_{6}$), or
\item
$N_{G}(a)\cap V(H)=\{u_{4},u_{5}\}$ (and so $G[V(H)\cup \{a\}]\in \H_{5}$).
\end{enumerate}
\end{lem}
\proof
For each $i\in \{2,3\}$, we note that the graph $B_{i}:=H-u_{i}$ belongs to $\H_{1}$.
Since $N_{G}(a)\cap V(H)\not=\emptyset $, $N_{G}(a)\cap V(B_{i})\not=\emptyset $ for some $i\in \{2,3\}$.
If $au_{4},au_{5}\not\in E(G)$, then $N_{G}(a)\cap V(B_{i})\subseteq \{u_{1},u_{u_{5-i}}\}\cup C_{u_{6}}$ for each $i\in \{2,3\}$, which contradicts Lemma~\ref{lem2.2}.
Thus, $au_{4}\in E(G)$ or $au_{5}\in E(G)$.
We may assume that $au_{4}\in E(G)$.
Then by Lemma~\ref{lem2.2}, we have either $G[V(B_{3})\cup \{a\}]\in \H_{1}$, or $N_{G}(a)\cap V(B_{3})=\{u_{1},u_{4}\}$ and $|C_{u_{6}}|=1$, or $N_{G}(a)\cap V(B_{3})=\{u_{4},u_{5}\}$.

\medskip
\noindent
\textbf{Case 1:} $G[V(B_{3})\cup \{a\}]\in \H_{1}$.

In this case, we have $\{u_{2},u_{4}\}\cup C_{u_{6}}\subseteq N_{G}(a)$.
Then again by Lemma~\ref{lem2.2}, we have either $G[V(B_{2})\cup \{a\}]\in \H_{1}$ or $N_{G}(a)\cap V(B_{2})=\{u_{1},u_{4}\}\cup C_{u_{6}}$ or $N_{G}(a)\cap V(B_{2})=\{u_{1},u_{3},u_{4}\}\cup C_{u_{6}}$.
If $N_{G}(a)\cap V(B_{2})=\{u_{1},u_{4}\}\cup C_{u_{6}}$ (i.e., $N_{G}(a)\cap V(H)=\{u_{1},u_{2},u_{4}\}\cup C_{u_{6}}$), then $G[\{u_{2},a,u_{1},u_{3},u_{5}\}]\cong Z_{2}$; if $N_{G}(a)\cap V(B_{2})=\{u_{1},u_{3},u_{4}\}\cup C_{u_{6}}$ (i.e., $N_{G}(a)\cap V(H)=\{u_{1},u_{2},u_{3},u_{4}\}\cup C_{u_{6}}$), then $G[\{u_{2},u_{4},a,u_{3},u_{5}\}]\cong Z_{2}$.
In either case, we get a contradiction.
Thus $G[V(B_{2})\cup \{a\}]\in \H_{1}$.
Since $au_{4}\in E(G)$, we see that $N_{G}(a)\cap V(B_{2})=\{u_{3},u_{4},u_{5}\}\cup C_{u_{6}}$, and hence $N_{G}(a)\cap V(H)=\{u_{2},u_{3},u_{4},u_{5}\}\cup C_{u_{6}}$.
Consequently, $G[V(H)\cup \{a\}]=H_{3}(\{C'_{u_{6}}\})$ where $C'_{u_{6}}=C_{u_{6}}\cup \{a\}$.

\medskip
\noindent
\textbf{Case 2:} $N_{G}(a)\cap V(B_{3})=\{u_{1},u_{4}\}$ and $|C_{u_{6}}|=1$.

Since $au_{2}\not\in E(G)$ and $G[\{u_{1},u_{2},u_{3},a\}]\not\cong K_{1,3}$, we have $au_{3}\in E(G)$.
Hence $N_{G}(a)\cap V(H)=\{u_{1},u_{3},u_{4}\}$ and $|C_{u_{6}}|=1$.

\medskip
\noindent
\textbf{Case 3:} $N_{G}(a)\cap V(B_{3})=\{u_{4},u_{5}\}$.

Since $G[\{u_{3},u_{1},b,a\}]  \not\cong K_{1,3}$ for any $b\in C_{u_{6}}$, we have $au_{3}\not\in E(G)$.
Hence $N_{G}(a)\cap V(H)=\{u_{4},u_{5}\}$.
\qed

\begin{lem}
\label{lem2.5}
Let $G$ be a connected $\{K_{1,3},Z_{2},N\}$-free graph which contains an induced subgraph $H=H_{4}(\{C_{v_{4}},C_{v_{5}},C_{v_{6}}\})$, where $C_{v_{4}}$, $C_{v_{5}}$ and $C_{v_{6}}$ are vertex-disjoint cliques.
Then for each vertex $a\in V(G)-V(H)$ with $N_{G}(a)\cap V(H)\not=\emptyset $, $G[V(H)\cup \{a\}]\in \H_{4}$.
\end{lem}
\proof
For each $i\in \{4,5,6\}$, let $b_{i}\in C_{v_{i}}$.
For each $i\in \{5,6\}$, we note that the graph $B_{i}:=H-C_{v_{i}}$ belongs to $\H_{1}$.

Suppose $N_{G}(a)\cap \{v_{1},v_{2},v_{3}\}=\emptyset $.
Since $N_{G}(a)\cap V(H)\not=\emptyset $, we may assume that $ab_{4}\in E(G)$.
Then $G[\{b_{4},a,v_{1},v_{2}\}] \cong K_{1,3}$, giving a contradiction.
Thus, $N_{G}(a)\cap \{v_{1},v_{2},v_{3}\}\not=\emptyset $.
We may assume that $av_{1}\in E(G)$.
Then, by Lemma~\ref{lem2.2}, we have $G[V(B_{5})\cup \{a\}]\in \H_{1}$ or $N_{G}(a)\cap V(B_{5})=\{v_{1},v_{i}\}$ for some $i\in\{2,3\}$.

Suppose that $N_{G}(a)\cap V(B_{5})=\{v_{1},v_{i}\}$ for some $i\in\{2,3\}$.
In this case, we may assume that $N_{G}(a)\cap V(B_{5})=\{v_{1},v_{2}\}$.
Then by Lemma~\ref{lem2.2}, $N_{G}(a)\cap V(B_{6})=\{v_{1},v_{2}\}$.
In particular, $N_{G}(a)\cap V(H)=\{v_{1},v_{2}\}$.
Then $G[\{b_{5},v_{3},b_{6},v_{1},a\}] \cong Z_{2}$, giving a contradiction.
Thus $G[V(H)\cup \{a\}]\in \H_{1}$.

Hence we have $N_{G}(a)\cap V(B_{5})=\{v_{1}\}\cup C_{v_{4}}\cup C_{v_{6}}$ or $N_{G}(a)\cap V(B_{5})=\{v_{1},v_{i}\}\cup C_{v_{4}}\cup C_{v_{6}}$ for some $i\in \{2,3\}$.
If $N_{G}(a)\cap V(B_{5})=\{v_{1}\}\cup C_{v_{4}}\cup C_{v_{6}}$, then $ab_{5}\in E(G)$ because $G[\{a,v_{1},b_{6},b_{5},v_{2}\}] \not\cong Z_{2}$, and hence $G[\{b_{5},a,v_{2},v_{3}\}] \cong K_{1,3}$, giving a contradiction.
Thus $N_{G}(a)\cap V(B_{5})=\{v_{1},v_{i}\}\cup C_{v_{4}}\cup C_{v_{6}}$ for some $i\in \{2,3\}$.
We may assume that $N_{G}(a)\cap V(B_{5})=\{v_{1},v_{2}\}\cup C_{v_{4}}\cup C_{v_{6}}$.
Since $\{v_{1},v_{2}\}\cup C_{v_{4}}\subseteq N_{G}(a)\cap V(B_{6})\subseteq \{v_{1},v_{2}\}\cup C_{v_{4}}\cup C_{v_{5}}$, we have $C_{v_{5}}\subseteq N_{G}(a)$ by Lemma~\ref{lem2.2}.
In particular, $N_{G}(a)\cap V(H)=\{v_{1},v_{2}\}\cup C_{v_{4}}\cup C_{v_{5}}\cup C_{v_{6}}$.
Therefore $G[V(H)\cup \{a\}]=H_{4}(\{C'_{v_{4}},C'_{v_{5}},C'_{v_{6}}\})$ where $C'_{s_{4}}=C_{s_{4}}\cup \{a\}$ and $C'_{s_{i}}=C_{s_{i}}~(i\in \{5,6\})$.
\qed

\begin{lem}
\label{lem2.6}
Let $G$ be a connected $\{K_{1,3},Z_{2},N\}$-free graph which contains an induced subgraph $H=H_{5}(\{C_{w_{7}}\})$, where $C_{w_{7}}$ is a clique.
Then for each vertex $a\in V(G)-V(H)$ with $N_{G}(a)\cap V(H)\not=\emptyset $, one of the following holds:
\begin{enumerate}[{\upshape(i)}]
\item
$G[V(H)\cup \{a\}]\in \H_{5}$, or
\item
$N_{G}(a)\cap V(H)=\{w_{1},w_{2},w_{i},w_{9-i}\}$ for some $i\in \{3,4\}$ and $|C_{w_{7}}|=1$ $($and so $G[V(H)\cup \{a\}]\cong H_{7})$.
\end{enumerate}
\end{lem}
\proof
For each $i\in \{1,2\}$, we note that the graph $B_{i}:=H-w_{i}$ belongs to $\H_{3}$.
Since $N_{G}(a)\cap V(H)\not=\emptyset $, $N_{G}(a)\cap V(B_{i})\not=\emptyset $ for some $i\in \{1,2\}$.
We may assume that $N_{G}(a)\cap V(B_{1})\not=\emptyset $.
If $N_{G}(a)\cap V(B_{1})=\{w_{3},w_{5}\}$, then either $N_{G}(a)\cap V(B_{2})=\{w_{3},w_{5}\}$ or $N_{G}(a)\cap V(B_{2})=\{w_{1},w_{3},w_{5}\}$, which contradicts Lemma~\ref{lem2.4}.
This, together with Lemma~\ref{lem2.4}, implies that either $G[V(B_{1})\cup \{a\}]\in \H_{3}$ or $N_{G}(a)\cap V(B_{1})=\{w_{2},w_{i},w_{9-i}\}$ for some $i\in \{3,4\}$ and $|C_{w_{7}}|=1$.

\medskip
\noindent
\textbf{Case 1:} $G[V(B_{1})\cup \{a\}]\in \H_{3}$.

Note that we have either $N_{G}(a)\cap V(B_{2})=\{w_{3},w_{4},w_{5},w_{6}\}\cup C_{w_{7}}$ or $N_{G}(a)\cap V(B_{2})=\{w_{1},w_{3},w_{4},w_{5},w_{6}\}\cup C_{w_{7}}$.
This, together with Lemma~\ref{lem2.4}, leads to $N_{G}(a)\cap V(H)=\{w_{3},w_{4},w_{5},w_{6}\}\cup C_{w_{7}}$.
Hence, $G[V(H)\cup \{a\}]=H_{5}(\{C'_{w_{7}}\})\in \H_{5}$ where $C'_{w_{7}}=C_{w_{7}}\cup \{a\}$.

\medskip
\noindent
\textbf{Case 2:} $N_{G}(a)\cap V(B_{1})=\{w_{2},w_{i},w_{9-i}\}$ for some $i\in \{3,4\}$ and $|C_{w_{7}}|=1$.

We may assume that $N_{G}(a)\cap V(B_{1})=\{w_{2},w_{3},w_{6}\}$.
Then $N_{G}(a)\cap V(B_{2})=\{w_{3},w_{6}\}$ or $N_{G}(a)\cap V(B_{2})=\{w_{1},w_{3},w_{6}\}$.
This, together with Lemma~\ref{lem2.4}, leads to $N_{G}(a)\cap V(H)=\{w_{1},w_{2},w_{3},w_{6}\}$.
\qed

\begin{lem}
\label{lem2.7}
Let $G$ be a connected $\{K_{1,3},Z_{2},N\}$-free graph which contains an induced subgraph $H=H_{6}$.
Then for each vertex $a\in V(G)-V(H)$ with $N_{G}(a)\cap V(H)\not=\emptyset $, $N_{G}(a)\cap V(H)=\{x_{i},x_{i+1},x_{7}\}$ for some $i\in \{1,3\}$. 
Consequently, $G[V(H)\cup \{a\}]\cong H_{7}$.
\end{lem}
\proof
We note that the graph $B:=H-x_{1}$ belongs to $\H_{3}$, and the graph $B^{*}:=H-x_{5}$ belongs to $\H_{2}$.

We first suppose that $ax_{i},ax_{i+2}\in E(G)$ for some $i\in \{1,2\}$.
We may assume that $ax_{1},ax_{3}\in E(G)$.
Then by Lemma~\ref{lem2.3}, we have $N_{G}(a)\cap V(B^{*})=\{x_{1},x_{3},x_{6},x_{7}\}$, and hence either $N_{G}(a)\cap V(B)=\{x_{3},x_{6},x_{7}\}$ or $N_{G}(a)\cap V(B)=\{x_{3},x_{5},x_{6},x_{7}\}$, which contradicts Lemma~\ref{lem2.4}.
Thus,
\begin{align}
\label{lem2.7 cond}
\mbox{for each $i\in \{1,2\}$, either $ax_{i}\not\in E(G)$ or $ax_{i+2}\not\in E(G)$.}
\end{align}

If $N_{G}(a)\cap V(B^{*})\not=\emptyset $, then $|N_{G}(a)\cap V(B^{*})|\geq 2$ by Lemma~\ref{lem2.3}.
In particular, we have $N_{G}(a)\cap V(B)\not=\emptyset $.
If $N_{G}(a)\cap \{x_{1},x_{2},x_{3},x_{4}\}=\emptyset $, then $N_{G}(a)\cap V(B)\subseteq \{x_{5},x_{6},x_{7}\}$, which contradicts Lemma~\ref{lem2.4}.
Thus we may assume that $ax_{1}\in E(G)$.
By (\ref{lem2.7 cond}), $ax_{3}\not\in E(G)$.
Since $G[\{x_{1},a,x_{2},x_{3}\}]\not\cong K_{1,3}$, we have $ax_{2}\in E(G)$. 
So, $ax_{4}\not\in E(G)$ by (\ref{lem2.7 cond}).
Then, by Lemma~\ref{lem2.4}, $N_{G}(a)\cap V(B)=\{x_{2},x_{7}\}$.
Consequently, $N_{G}(a)\cap V(H)=\{x_{1},x_{2},x_{7}\}$.
\qed

\begin{lem}
\label{lem2.8}
Let $G$ be a connected $\{K_{1,3},Z_{2},N\}$-free graph which contains an induced subgraph $H=H_{7}$.
Then for each vertex $a\in V(G)-V(H)$ with $N_{G}(a)\cap V(H)\not=\emptyset $, $N_{G}(a)\cap V(H)=\{y_{1},y_{2},y_{7},y_{8}\}$.
Consequently, $G[V(H)\cup \{a\}]\cong H_{8}$.
\end{lem}
\proof
For each $i\in \{1,2\}$, we note that the graph $B_{i}:=H-y_{i}$ is isomorphic to $H_{6}$.
Since $N_{G}(a)\cap V(H)\not=\emptyset $, we have $N_{G}(a)\cap V(B_{i})\not=\emptyset $ for some $i\in \{1,2\}$.
We may assume that $N_{G}(a)\cap V(B_{1})\not=\emptyset $.
Then, by Lemma~\ref{lem2.7}, $N_{G}(a)\cap V(B_{1})=\{y_{2},y_{3},y_{5}\}$ or $N_{G}(a)\cap V(B_{1})=\{y_{2},y_{7},y_{8}\}$.
In particular, $\{y_{3},y_{5}\}\subseteq N_{G}(a)\cap V(B_{2})$ or $\{y_{7},y_{8}\}\subseteq N_{G}(a)\cap V(B_{2})$.
This, together with Lemma~\ref{lem2.7}, leads to $N_{G}(a)\cap V(B_{1})=\{y_{2},y_{7},y_{8}\}$ and $N_{G}(a)\cap V(B_{2})=\{y_{1},y_{7},y_{8}\}$. 
So, $N_{G}(a)=\{y_{1},y_{2},y_{7},y_{8}\}$.
\qed

\medbreak\noindent\textit{Proof of Theorem~\ref{thm1z}.}\quad
By routine but tedious argument, we can verify that every graph in $\bigcup _{0\leq i\leq 8}\H_{i}$ is $\{K_{1,3},Z_{2}\}$-free but not $B_{1,1}$-free (and we omit its detail).
Thus it suffices to show that, if a connected $\{K_{1,3},Z_{2}\}$-free graph $G$ is not $B_{1,1}$-free (i.e., $G$ contains $B_{1,1}$ as an induced subgraph), then $G$ belongs to $\bigcup _{0\leq i\leq 8}\H_{i}$.

Assume that $G$ contains $B_{1,1}~(\in \H_{1})$ as an induced subgraph.
Then $G$ contains a graph $H\in \bigcup _{0\leq i\leq 8}\H_{i}$ as an induced subgraph.
Choose $H$ so that $|V(H)|$ is as large as possible.
It suffices to show that $G=H$.
By way of contradiction, suppose that $G\not=H$ (i.e., $V(G)-V(H)\not=\emptyset $).
Then by Lemma~\ref{lem2.1}, $G$ is $N$-free.
Since $G$ is connected, there exists a vertex $a\in V(G)-V(H)$ which is adjacent to a vertex in $V(H)$.
By the maximality of $H$, $G[V(H)\cup \{a\}]\not\in \bigcup _{0\leq i\leq 8}\H_{i}$.
This, together with Lemmas~\ref{lem2.2}--\ref{lem2.8}, gives $H=H_{8}$.
For each $i\in \{7,8,9\}$, we note that the graph $B_{i}:=H-z_{i}$ is isomorphic to $H_{7}$.
Since $N_{G}(a)\cap V(H)\not=\emptyset $, we have $N_{G}(a)\cap V(B_{i})\not=\emptyset $ for some $i\in \{7,8\}$.
We may assume that $N_{G}(a)\cap V(B_{7})\not=\emptyset $.
Then by Lemma~\ref{lem2.7}, $N_{G}(a)\cap V(B_{7})=\{z_{3},z_{4},z_{8},z_{9}\}$.
In particular, $az_{3}\in E(G)$.
On the other hand, since $N_{G}(a)\cap V(B_{9})\not=\emptyset $, $N_{G}(a)\cap V(B_{9})=\{z_{1},z_{2},z_{7},z_{8}\}$, and so $az_{3}\not\in E(G)$, giving a contradiction.

This completes the proof of Theorem~\ref{thm1z}.
\qed

\subsection{Proof of Theorem~\ref{thm3}}\label{subsec2.2}

In order to prove Theorem~\ref{thm3}, we give a further definition.
For two integers $s$ and $t$, we let $[s,t]=\{i\in \mathbb{N}\mid s\leq i\leq t\}$.
Note that if $s>t$, then $[s,t]=\emptyset $.

Here we prove Theorem~\ref{thm3}.
We can easily verify that every graph in $\P(\max\{3m,m+4\})$ is $\{K_{1,3},B_{1,m}\}$-free but not $P_{\max\{3m,m+4\}}$-free.
Thus it suffices to show that if a connected $\{K_{1,3},B_{1,m}\}$-free graph $G$ is not $P_{\max\{3m,m+4\}}$-free (i.e., $G$ contains $P_{\max\{3m,m+4\}}$ as an induced subgraph), then $G$ belongs to $\P(\max\{3m,m+4\})$.

Assume that $G$ contains $P_{\max\{3m,m+4\}}$ as an induced subgraph.
Then $G$ contains a graph $H\in \P(\max\{3m,m+4\})$ as an induced subgraph.
Choose $H$ so that $|V(H)|$ is as large as possible.
It suffices to show that $G=H$.
Otherwise, there exists a vertex $a\in V(G)-V(H)$ such that $N_G(a) \cap V(H) \neq \emptyset$.
Let $l$ be the integer so that $H$ is either a fat $l$-path or a fat $l$-cycle.
Then we can write either $H=F_{p}(L_{1},\ldots ,L_{l})$ or $H=F_{c}(L_{0},\ldots ,L_{l})$ for some vertex-disjoint cliques $L_{0},\ldots ,L_{l}$.
Let $I=\{i\mid N_{G}(a)\cap L_{i}\not=\emptyset \}$.

\begin{claim}
\label{cl3.1}
$|I|\leq 4$.
\end{claim}
\proof
Suppose that there are five fundamental cliques $L^{(1)},\ldots ,L^{(5)}$ of $H$ with $N_{G}(a)\cap L^{(i)}\not=\emptyset ~(1\leq i\leq 5)$.
For each $i~(1\leq i\leq 5)$, let $b^{(i)}\in N_{G}(a)\cap L^{(i)}$.
Since $\max\{3m,m+4\}\geq 5$, if $H$ is a fat cycle, then $H$ has at least six fundamental cliques.
Thus $G[\{b^{(i)}\mid 1\leq i\leq 5\}]$ has no cycle, and so is a forest of order five and maximum degree at most two.
Then we can easily check that $G[\{b^{(i)}\mid 1\leq i\leq 5\}]$ has an independent set $B$ with $|B|=3$, and hence $G[\{a\}\cup B]\cong K_{1,3}$, giving a contradiction.
\qed

If $H$ is a fat cycle, then $N_{G}(a)\cap L_{i}=\emptyset $ for some $0\leq i\leq l$ by Claim~\ref{cl3.1}. 
By relabeling $L_{0},\ldots ,L_{l}$ if necessary, we may assume that
\begin{enumerate}[(L1)]
\item $0\not\in I$, and
\item subject to (L1), $|I\cap \{1,l\}|$ is as small as possible.
\end{enumerate}
Thus, if $H$ is a fat cycle and there exists an integer $i~(1\leq i\leq l-2)$ with $i,i+1,i+2\not\in I$, then $I\cap \{0,1,l\}=\emptyset $.

For each $i~(1\leq i\leq l)$, we take a vertex $b_{i}$ as follows:
If $i\in I$, let $b_{i}\in N_{G}(a)\cap L_{i}$; otherwise (i.e., $i\not\in I$), let $b_{i} \in L_{i}$.
Note that, by our choices of indices, $b_{1}b_{l}\not\in E(G)$ regardless of $H$ being a fat path or a fat cycle.

\begin{claim}
\label{cl3.2}
Assume that there exists an index $j~(2\leq j\leq l-2)$ such that $I\cap [2,l-1]=\{j,j+1\}$.
Then either $j=2$ and $ab_{1}\in E(G)$ or $j=l-2$ and $ab_{l}\in E(G)$.
\end{claim}
\proof
Recall that $l \ge \max\{ 3m, m+4\}$.
We first consider the case $l-m-1\leq j\leq m+1$.
Then $l\leq 2m+2$.
Since $l\geq 3m$, we have $m\leq 2$; since $l\geq m+4$, we have $m\geq 2$.
Hence $m=2$, and this forces $l=6$ and $j=3$.
By the assumption of the claim, $ab_{2},ab_{5}\not\in E(G)$.
Since $G[\{b_{2},b_{3},a,b_{4},b_{5},b_{6}\}]\not\cong B_{1,2}$ and $G[\{b_{5},b_{4},a,b_{3},b_{2},b_{1}\}]\not\cong B_{1,2}$, 
we have $ab_{1},ab_{6}\in E(G)$.
Then $G[\{a,b_{1},b_{3},b_{6}\}] \cong K_{1,3}$, giving a contradiction.
Thus either $j\geq m+2$ or $j\leq l-m-2$.

We now consider the case $j\geq m+2$ (i.e., $j-m\geq 2$).
Then $ab_{i}\not\in E(G)$ for every $j-m\leq i\leq j-1$.
Since $G[\{b_{j+2},b_{j+1},a,b_{j},b_{j-1},\ldots ,b_{j-m}\}]\not\cong B_{1,m}$, this forces $b_{j+2}=b_{l}$ (i.e., $j=l-2$) and $ab_{l}\in E(G)$, as desired.
Thus we may assume that $j\leq l-m-2$ (i.e., $j+m+1\leq l-1$).
Then $ab_{i}\not\in E(G)$ for every $j+2\leq i\leq j+m+1$.
Since $G[\{b_{j-1},b_{j},a,b_{j+1},b_{j+2},\ldots ,b_{j+m+1}\}]\not\cong B_{1,m}$, this forces $b_{j-1}=b_{1}$ (i.e., $j=2$) and $ab_{1}\in E(G)$, as desired.
\qed

\begin{claim}
\label{cl3.3}
For each $j\in I$, 
there exists an index $j'~(j'\not=j)$ such that $|j-j'|=1$ and $j'\in I$.
\end{claim}
\proof
If $2\leq j\leq l-1$ and $j-1,j+1\not\in I$, then $G[\{b_{j},a,b_{j-1},b_{j+1}\}] \cong K_{1,3}$, giving a contradiction.
Hence if $2\leq j\leq l-1$, then the desired conclusion holds.
Thus we may assume that $j\in \{1,l\}$.

For the moment, we assume that $j=1$ and $2\not\in I$.
We further suppose that there exists an index $i~(3\leq i\leq l-1)$ with $i\in I$.
Choose $i$ so that $i$ is as small as possible.
Then, $i+1\in I$ since $3\leq i\leq l-1$.
If $i+1\leq l-1$ and $I\cap [3,l-1]=\{i,i+1\}$, then $i=l-2$ and $l\in I$ by Claim~\ref{cl3.2}.
This implies that if $i+1\leq l-1$ (i.e., $i\leq l-2$), then there are three indices $i_{1},i_{2},i_{3}~(3\leq i_{1}<i_{2}<i_{3}\leq l)$ with $i_{1},i_{2},i_{3}\in I$, and hence $G[\{a,b_{1},b_{i_{1}},b_{i_{3}}\}] \cong K_{1,3}$, giving a contradiction.
Thus $i\geq l-1$, and so $i=l-1$.
Note that $I\cap [1,l]=\{1,l-1,l\}$.
This, together with the fact $l-m-1\geq 3$, implies that $G[\{b_{1},a,b_{l},b_{l-1},\ldots ,b_{l-m-1}\}] \cong B_{1,m}$, giving a contradiction.
Thus $I\cap [2,l-1]=\emptyset $.
By the choice of $L_{i}$, we see that $H$ is a fat path.
If there exists a vertex $u\in L_{1}$ with $au\not\in E(G)$, then, since $l\geq m+4$, $G[\{a,b_{1},u,b_{2},\ldots ,b_{m+2}\}] \cong B_{1,m}$, giving a contradiction.
Thus $L_{1}\subseteq N_{G}(a)$.
By the symmetry and the fact $l-1\not\in I$, if $l\in I$, then $L_{l}\subseteq N_{G}(a)$.
Hence either $N_{G}(a)\cap V(H)=L_{1}$ or $N_{G}(a)\cap V(H)=L_{1}\cup L_{l}$.
If $N_{G}(a)\cap V(H)=L_{1}$, then $G[V(H)\cup \{a\}]=F_{p}(\{a\},L_{1},\ldots ,L_{l})$; if $N_{G}(a)\cap V(H)=L_{1}\cup L_{l}$, then $G[V(H)\cup \{a\}]=F_{c}(\{a\},L_{1},\ldots ,L_{l})$.
In either case, $G[V(H)\cup \{a\}]\in \P(\max\{3m,m+4\})$, which contradicts the maximality of $H$.
Thus if $j=1$, then $2\in I$.
By the symmetry, if $j=l$, then $l-1\in I$.
\qed

Let $i_{1}=\min \{i\mid i\in I\}$ and $i_{2}=\max \{i\mid i\in I\}$.
By Claim~\ref{cl3.3}, $i_{1}+1,i_{2}-1\in I$.

\begin{claim}
\label{cl3.4}
If $i_{1}\not=1$, then $L_{i_{1}+1}\subseteq N_{G}(a)$.
If $i_{2}\not=l$, then $L_{i_{2}-1}\subseteq N_{G}(a)$.
\end{claim}
\proof
If $i_{1}\not=1$ and $L_{i_{1}+1}\not\subseteq N_{G}(a)$, say $u\in L_{i_{1}+1}-N_{G}(a)$, then $G[\{b_{i_{1}},b_{i_{1}-1},a,u\}]\cong K_{1,3}$, giving a contradiction.
Thus if $i_{1}\not=1$, then $L_{i_{1}+1}\subseteq N_{G}(a)$.
By the symmetry, we have $L_{i_{2}-1}\subseteq N_{G}(a)$ if $i_{2}\not=l$.
\qed

Since $|I|\leq 4$,
we divide the rest of the proof into three cases 
according to $|I| \le 2$, $|I|=3$, and $|I|=4$.

\medskip
\noindent
\textbf{Case 1:} $|I|\leq 2$.

By Claim~\ref{cl3.3}, $I=\{i_{1},i_{2}\}=\{i_{1},i_{1}+1\}$.
If $|I\cap [2,l-1]|=2$, then either $1\in I$ or $l\in I$ by Claim~\ref{cl3.2}, and so $|I|\geq 3$, giving a contradiction.
Thus $|I\cap [2,l-1]|\leq 1$, which implies either $I=\{1,2\}$ or $I=\{l-1,l\}$.
We may assume that $I=\{1,2\}$.
By the choice of $L_{i}$, $H$ is a fat path.
If $L_{2}\not\subseteq N_{G}(a)$, say $u\in L_{2}-N_{G}(a)$, then $G[\{a,b_{2},u,b_{3},b_{4},\ldots ,b_{m+3}\}] \cong B_{1,m}$, giving a contradiction.
Thus $L_{2}\subseteq N_{G}(a)$.
This, together with Claim~\ref{cl3.4}, leads to $N_{G}(a)\cap V(H)=L_{1}\cup L_{2}$, and hence $G[V(H)\cup \{a\}]=F_{p}(L_{1}\cup \{a\},L_{2},\ldots ,L_{l})\in \P(\max\{3m,m+4\})$, which contradicts the maximality of $H$.

\medskip
\noindent
\textbf{Case 2:} $|I|=3$.

In this case, $I=\{i_{1},i_{1}+1~(=i_{2}-1),i_{2}\}=\{i_{1},i_{1}+1,i_{1}+2\}$.
By the choice of $L_{i}$, $H$ is a fat path.
Since either $i_{1}\not=1$ or $i_{2}\not=l$, $L_{i_{1}+1}~(=L_{i_{2}-1})\subseteq N_{G}(a)$ by Claim~\ref{cl3.4}.
Suppose that either $L_{i_{1}}\not\subseteq N_{G}(a)$ or $L_{i_{2}}\not\subseteq N_{G}(a)$.
We may assume that $L_{i_{1}}\not\subseteq N_{G}(a)$.
Let $u\in L_{i_{1}}-N_{G}(a)$.
If $i_{1}\leq l-m-2$ (i.e., $i_{1}+m+2\leq l$), then $G[\{u,b_{i_{1}+1},a,b_{i_{1}+2},\ldots ,b_{i_{1}+m+2}\}] \cong B_{1,m}$; if $i_{1}\geq m+2$ (i.e., $i_{1}-m-1\geq 1$), then $G[\{a,b_{i_{1}},u,b_{i_{1}-1},\ldots ,b_{i_{1}-m-1}\}] \cong B_{1,m}$.
In either case, we get a contradiction.
Thus $l-m-1\leq i_{1}\leq m+1$.
This, together with the assumption $l\geq \max\{3m,m+4\}$, leads to $m=2$, $l=6$ and $i_{1}=3$.
Then $G[\{b_{1},b_{2},u,b_{3},a,b_{5}\}]\cong B_{1,2}$, giving a contradiction.
Thus $L_{i_{1}}\cup L_{i_{2}}\subseteq N_{G}(a)$ (i.e., $N_{G}(a)\cap V(H)=L_{i_{1}}\cup L_{i_{1}+1}\cup L_{i_{2}}$).
Hence $G[V(H)\cup \{a\}]=F_{p}(L_{1},\ldots ,L_{i_{1}},L_{i_{1}+1}\cup \{a\},L_{i_{2}},\ldots ,L_{l})$, which contradicts the maximality of $H$.

\medskip
\noindent
\textbf{Case 3:} $|I|=4$.

In this case, $i_{1}+1<i_{2}-1$ and $I=\{i_{1},i_{1}+1,i_{2}-1,i_{2}\}$.
Let $J_{1}=[1,i_{1}-1]$, $J_{2}=[i_{1}+2,i_{2}-2]$ and $J_{3}=[i_{2}+1,l]$ (where $J_{i}$ may be empty).
If $|J_{1}|\geq m$, then $i_{1}-m\geq 1$, and hence $G[\{b_{i_{2}},a,b_{i_{1}+1},b_{i_{1}},b_{i_{1}-1},\ldots ,b_{i_{1}-m}\}] \cong B_{1,m}$; if $|J_{2}|\geq m$, then $i_{1}+m+1\leq i_{2}-2$, and hence $G[\{b_{i_{2}},a,b_{i_{1}},b_{i_{1}+1},\ldots ,b_{i_{1}+m+1}\}] \cong B_{1,m}$; if $|J_{3}|\geq m$, then $i_{2}+m\leq l$, and hence $G[\{b_{i_{1}},a,b_{i_{2}-1},b_{i_{2}},b_{i_{2}+1},\ldots ,b_{i_{2}+m}\}] \cong B_{1,m}$.
In either case, we get a contradiction.
Thus $\max \{|J_{1}|,|J_{2}|,|J_{3}|\}\leq m-1$.
On the other hand, $|J_{1}|+|J_{2}|+|J_{3}|=|[1,l]-\{i_{1},i_{1}+1,i_{2}-1,i_{2}\}|=l-4\geq \max\{3m-4,m\}$.
Hence we see that $m\geq 2$ and $|J_{i}|=|J_{i'}|=m-1$ for some $i,i'\in \{1,2,3\}$ with $i\not=i'$.
Without loss of generality, we may assume that $|J_{1}|=m-1$ (i.e., $i_{1}=m$).
If $|J_{2}|=m-1$, then $G[\{b_{i_{2}-2},b_{i_{2}-1},b_{i_{2}},a,b_{m},b_{m-1},\ldots ,b_{1}\}] \cong B_{1,m}$; if $|J_{3}|=m-1$, then $G[\{b_{i_{2}+1},b_{i_{2}},b_{i_{2}-1},a,b_{m},b_{m-1},\ldots ,b_{1}\}] \cong B_{1,m}$.
In either case, we again get a contradiction.

This completes the proof of Theorem~\ref{thm3}.

\section*{Appendix}

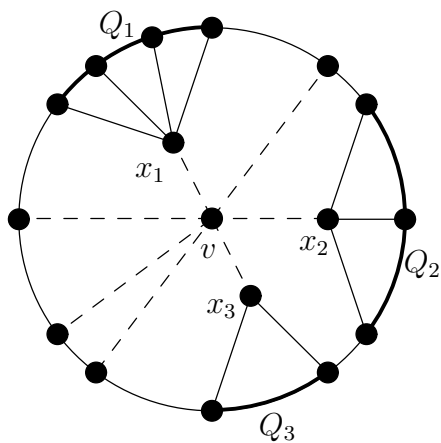
\begin{figure}[htb]
\begin{center}
{\unitlength 0.1in
\begin{picture}( 21.0800, 21.3300)( 11.4600,-28.5800)
%
\special{sh 1.000}%
\special{ia 2200 804 54 54  0.0000000  6.2831853}%
\special{pn 8}%
\special{ar 2200 804 54 54  0.0000000  6.2831853}%
%
\special{pn 8}%
\special{ar 2200 1804 1000 1000  0.0000000  6.2831853}%
%
\special{sh 1.000}%
\special{ia 1200 1804 54 54  0.0000000  6.2831853}%
\special{pn 8}%
\special{ar 1200 1804 54 54  0.0000000  6.2831853}%
%
\special{sh 1.000}%
\special{ia 3200 1804 54 54  0.0000000  6.2831853}%
\special{pn 8}%
\special{ar 3200 1804 54 54  0.0000000  6.2831853}%
%
\special{sh 1.000}%
\special{ia 2200 2804 54 54  0.0000000  6.2831853}%
\special{pn 8}%
\special{ar 2200 2804 54 54  0.0000000  6.2831853}%
%
\special{sh 1.000}%
\special{ia 2800 2604 54 54  0.0000000  6.2831853}%
\special{pn 8}%
\special{ar 2800 2604 54 54  0.0000000  6.2831853}%
%
\special{sh 1.000}%
\special{ia 2800 1004 54 54  0.0000000  6.2831853}%
\special{pn 8}%
\special{ar 2800 1004 54 54  0.0000000  6.2831853}%
%
\special{sh 1.000}%
\special{ia 3000 1204 54 54  0.0000000  6.2831853}%
\special{pn 8}%
\special{ar 3000 1204 54 54  0.0000000  6.2831853}%
%
\special{sh 1.000}%
\special{ia 3000 2404 54 54  0.0000000  6.2831853}%
\special{pn 8}%
\special{ar 3000 2404 54 54  0.0000000  6.2831853}%
%
\special{sh 1.000}%
\special{ia 1600 2604 54 54  0.0000000  6.2831853}%
\special{pn 8}%
\special{ar 1600 2604 54 54  0.0000000  6.2831853}%
%
\special{sh 1.000}%
\special{ia 1400 2404 54 54  0.0000000  6.2831853}%
\special{pn 8}%
\special{ar 1400 2404 54 54  0.0000000  6.2831853}%
%
\special{sh 1.000}%
\special{ia 1400 1204 54 54  0.0000000  6.2831853}%
\special{pn 8}%
\special{ar 1400 1204 54 54  0.0000000  6.2831853}%
%
\special{sh 1.000}%
\special{ia 1600 1004 54 54  0.0000000  6.2831853}%
\special{pn 8}%
\special{ar 1600 1004 54 54  0.0000000  6.2831853}%
%
\special{sh 1.000}%
\special{ia 2000 1404 54 54  0.0000000  6.2831853}%
\special{pn 8}%
\special{ar 2000 1404 54 54  0.0000000  6.2831853}%
%
\special{sh 1.000}%
\special{ia 2800 1804 54 54  0.0000000  6.2831853}%
\special{pn 8}%
\special{ar 2800 1804 54 54  0.0000000  6.2831853}%
%
\special{sh 1.000}%
\special{ia 2400 2204 54 54  0.0000000  6.2831853}%
\special{pn 8}%
\special{ar 2400 2204 54 54  0.0000000  6.2831853}%
%
\special{pn 8}%
\special{pa 2000 1404}%
\special{pa 1400 1204}%
\special{fp}%
\special{pa 2000 1404}%
\special{pa 1600 1004}%
\special{fp}%
\special{pa 2000 1404}%
\special{pa 2200 804}%
\special{fp}%
\special{pa 2800 1804}%
\special{pa 3000 1204}%
\special{fp}%
\special{pa 2800 1804}%
\special{pa 3200 1804}%
\special{fp}%
\special{pa 2800 1804}%
\special{pa 3000 2404}%
\special{fp}%
\special{pa 2800 2604}%
\special{pa 2400 2204}%
\special{fp}%
\special{pa 2400 2204}%
\special{pa 2200 2804}%
\special{fp}%
\put(18.8000,-15.6000){\makebox(0,0){$x_{1}$}}%
\put(27.3000,-19.4000){\makebox(0,0){$x_{2}$}}%
\put(22.5000,-22.7000){\makebox(0,0){$x_{3}$}}%
%
\special{sh 1.000}%
\special{ia 1890 854 54 54  0.0000000  6.2831853}%
\special{pn 8}%
\special{ar 1890 854 54 54  0.0000000  6.2831853}%
%
\special{pn 8}%
\special{pa 1890 854}%
\special{pa 2000 1404}%
\special{fp}%
\put(17.1000,-7.9000){\makebox(0,0){$Q_{1}$}}%
\put(32.9000,-20.4000){\makebox(0,0){$Q_{2}$}}%
\put(25.4000,-28.9000){\makebox(0,0){$Q_{3}$}}%
%
\special{pn 20}%
\special{ar 2200 1800 1000 1000  3.7850938  4.7123890}%
%
\special{pn 20}%
\special{ar 2200 1800 1000 1000  5.6396842  0.6435011}%
%
\special{pn 20}%
\special{ar 2200 1800 1000 1000  0.9272952  1.5707963}%
%
\special{sh 1.000}%
\special{ia 2200 1800 54 54  0.0000000  6.2831853}%
\special{pn 8}%
\special{ar 2200 1800 54 54  0.0000000  6.2831853}%
%
\special{pn 8}%
\special{pa 2200 1800}%
\special{pa 1200 1800}%
\special{da 0.070}%
\special{pa 2200 1800}%
\special{pa 1400 2400}%
\special{da 0.070}%
%
\special{pn 8}%
\special{pa 1600 2600}%
\special{pa 2200 1800}%
\special{da 0.070}%
\special{pa 2200 1800}%
\special{pa 2800 1000}%
\special{da 0.070}%
\special{pa 2800 1800}%
\special{pa 2200 1800}%
\special{da 0.070}%
\special{pa 2200 1800}%
\special{pa 2400 2200}%
\special{da 0.070}%
\special{pa 2200 1800}%
\special{pa 2000 1400}%
\special{da 0.070}%
\put(21.8000,-19.7000){\makebox(0,0){$v$}}%
\end{picture}}%

\caption{A fan-cycle system}
\label{fig_fc}
\end{center}
\end{figure}

Let $G$ be a graph.
A sequence $(C;v;Q_{1},\ldots ,Q_{m};x_{1},\ldots ,x_{m})$ is a {\it fan-cycle system} of $G$ if
\begin{enumerate}
\item
$C$ is a cycle of $G$,
\item
$Q_{1},\ldots ,Q_{m}$ are vertex disjoint paths of order at least two on $C$,
\item
$v,x_{1},\ldots ,x_{m}$ are distinct vertices with $V(G)-V(C)=\{v,x_{1},\ldots ,x_{m}\}$,
\item
$|V(C)-\bigcup _{1\leq i\leq m}V(Q_{i})|+m\geq 3$,
\item
$v$ is adjacent to every vertex in $(V(C)-\bigcup _{1\leq i\leq m}V(Q_{i}))\cup \{x_{1},\ldots ,x_{m}\}$, and
\item
for $i~(1\leq i\leq m)$, $x_{i}$ is adjacent to every vertex of $Q_{i}$
\end{enumerate}
(see Figure~\ref{fig_fc}).
In~\cite{CHOST}, the following lemma was proved in order to construct a spanning Halin subgraph.

\begin{lem}[Chen, Han, O, Shan and Tsuchiya~\cite{CHOST}]
\label{lm:fc}
If a graph $G$ has a fan-cycle system, then $G$ has a spanning Halin subgraph.
\end{lem}

Now we show that all $3$-connected graphs in $\bigcup_{0\le i\le 8}\H_{i}$ (defined in Subsection~\ref{subsec1.2}) have a spanning Halin subgraph.

\begin{lem}
\label{lm:appSHS}
For $G\in \bigcup_{0\le i\le 8}\H_{i}$, if $G$ is $3$-connected, then $G$ has a spanning Halin subgraph.
\end{lem}
\proof
Since all graphs in $\bigcup _{i\in \{2,3,5,6\}}\H_{i}$ are not $3$-connected, $G\in \H_{i}$ for some $i\in \{0,1,4,7,8\}$.
By Lemma~\ref{lm:fc}, it suffices to show that $G$ has a fan-cycle system.

\medskip
\noindent
\textbf{Case 1:} $G\in \H_{0}$.

Let $L_{1},\ldots ,L_{m}$ be the leaf-cliques of $G$, and let $R_{i}$ be the root of $L_{i}$.
For each $i~(1\leq i\leq m)$, let $v_{i}\in R_{i}$.
Since $G$ is $3$-connected, $|R_{i}-\{v_{i}\}|\geq 2$ for all $i$, and hence $G-\{v_{i}\mid 1\leq i\leq m\}$ has a Hamiltonian cycle $C$ containing $m-1$ vertex disjoint paths $Q_{2},\ldots ,Q_{m}$ with $V(Q_{i})=L_{i}\cup (R_{i}-\{v_{i}\})~(2\leq i\leq m)$.
Then $(C;v_{1};Q_{2},\ldots ,Q_{m};v_{2},\ldots ,v_{m})$ is a fan-cycle system of $G$.

\medskip
\noindent
\textbf{Case 2:} $G\in \H_{1}$.

Write $G=H_{1}(\{C_{s_{3}},C_{s_{4}},C_{s_{5}}\})$.
For each $i\in \{3,4\}$, let $a_{i}\in C_{s_{i}}$.
Since $G$ is $3$-connected, $|C_{s_{i}}| \geq 3$ for $i\in \{3,4\}$, and hence $G-\{a_{3},a_{4}\}$ has a Hamiltonian cycle $C$ containing a path $Q$ with $V(Q)=(C_{s_{4}}-\{a_{4}\})\cup \{s_{2}\}$.
Then $(C;a_{3};Q;a_{4})$ is a fan-cycle system of $G$.

\medskip
\noindent
\textbf{Case 3:} $G\in \H_{4}$.

Write $G=H_{1}(\{C_{v_{4}},C_{v_{5}},C_{v_{6}}\})$.
Since $G$ is $3$-connected, $|C_{v_{i}}\cup C_{v_{j}}|\geq 3$ for $i,j\in \{4,5,6\}$ with $i\not=j$.
By symmetry, we may assume that $|C_{v_{4}}|\geq 2$ and $|C_{v_{5}}|\geq 2$.
For each $i\in \{4,5\}$, let $a_{i}\in C_{v_{i}}$.
Then $G-\{a_{4},a_{5}\}$ has a Hamiltonian cycle $C$ containing a path $Q$ with $V(Q)=(C_{v_{5}}-\{a_{5}\})\cup \{v_{3}\}$, and hence $G$ has a fan-cycle system $(C;a_{4};Q;a_{5})$.

\medskip
\noindent
\textbf{Case 4:} $G\in \H_{7}$.

Let $C=y_{2}y_{4}y_{7}y_{8}y_{6}$ be a cycle of $G$, and $Q_{1}=y_{4}y_{7}$ and $Q_{2}=y_{8}y_{6}$ be paths on $C$.
Then $(C;y_{1};Q_{1},Q_{2};y_{3},y_{5})$ is a fan-cycle system of $G$.

\medskip
\noindent
\textbf{Case 5:} $G\in \H_{8}$.

Let $C=z_{2}z_{4}z_{7}z_{9}z_{8}z_{6}$ be a cycle of $G$, and $Q_{1}=z_{4}z_{7}$ and $Q_{2}=z_{8}z_{6}$ be paths on $C$.
Then $(C;z_{1};Q_{1},Q_{2};z_{3},z_{5})$ is a fan-cycle system of $G$.

This completes the proof of Lemma~\ref{lm:appSHS}.
\qed

\begin{lem}
\label{lm:appSHS2}
For $G\in \P(5)$, if $G$ is $3$-connected, then $G$ has a spanning Halin subgraph.
\end{lem}
\proof
We first suppose that $G$ is a fat path, and write $G=F_{p}(L_{1},\ldots ,L_{l})$.
For each $i~(2\leq i\leq l-1)$, let $a_{i}\in L_{i}$.
Since $G$ is $3$-connected, $|L_{i}-\{a_{i}\}|\geq 2$ for $i~(2\leq i\leq l-1)$, and hence $G-\{a_{2},\ldots ,a_{l-1}\}$ has a Hamiltonian cycle $C$ such that $C[L_{i}]$ has exactly two components for every $i~(2\leq i\leq l-1)$.
We take the spanning tree $T$ of $G$ such that $N_{T}(a_{2})=L_{1}\cup (L_{2}-\{a_{2}\})\cup \{a_{3}\}$, $N_{T}(a_{l-1})=L_{l}\cup (L_{l-1}-\{a_{l-1}\})\cup \{a_{l-2}\}$ and $N_{T}(a_{i})=(L_{i}-\{a_{i}\})\cup \{a_{i-1},a_{i+1}\}~(3\leq i\leq l-2)$.
Then $T$ has no vertices of degree $2$ and $V(G)-\{a_{2},\ldots ,a_{l-1}\}$ is the set of leaves of $T$.
Hence $T\cup C$ is a spanning Halin subgraph of $G$.

We next suppose that $G$ is a fat cycle, and write $G=F_{c}(L_{0},\ldots ,L_{l})$.
Since $G$ is $3$-connected, $G$ has at most two fundamental cliques of order one.
Furthermore, if $G$ has exactly two fundamental cliques of order one, then such cliques are consecutive.
By symmetry, we may assume that $|L_{i}|\geq 2$ for every $i~(1\leq i\leq l-1)$.
For each $i~(1\leq i\leq l-1)$, let $a_{i}\in L_{i}$.
Then $G-\{a_{1},\ldots ,a_{l-1}\}$ has a Hamiltonian cycle $C$ such that $C[L_{i}]$ has exactly one component for every $i~(0\leq i\leq l)$.
We take a spanning tree $T$ of $G$ such that $N_{T}(a_{1})=L_{0}\cup (L_{1}-\{a_{1}\})\cup \{a_{2}\}$, $N_{T}(a_{l-1})=L_{l}\cup (L_{l-1}-\{a_{l-1}\})\cup \{a_{l-2}\}$ and $N_{T}(a_{i})=(L_{i}-\{a_{i}\})\cup \{a_{i-1},a_{i+1}\}~(2\leq i\leq l-2)$.
Then $T$ has no vertices of degree $2$ and $V(G)-\{a_{1},\ldots ,a_{l-1}\}$ is the set of leaves of $T$.
Hence $T\cup C$ is a spanning Halin subgraph of $G$.
\qed

Theorems~\ref{thm1z}, \ref{thm:P5} and Lemma~\ref{lm:appSHS} lead to Theorem~\ref{thm:SHS2}.
Theorems~\ref{thm1}, \ref{thm:SHS2} and Lemma~\ref{lm:appSHS2} lead to Theorem~\ref{thm:SHS3}.


\begin{thebibliography}{99}
\bibitem {Be} 
P.~Bedrossian,
Forbidden subgraph and minimum degree conditions for hamiltonicity,
Ph.D.~Thesis, Memphis State University, 1991.

\bibitem{BH}
A.~Brandst\"{a}dt and P.L.~Hammer,
On the stability number of claw-free $P_{5}$-free and more general graphs,
Discrete Appl. Math. {\bf 95} (1999), 163--167. 

\bibitem{BV}
H.~Broersma and H.J.~Veldman,
Restrictions on induced subgraphs ensuring Hamiltonicity or pancyclicity of $K_{1,3}$-free graphs,
{\it Contemporary methods in graph theory}, Bibliographisches Inst., Mannheim (1990), 181--194.

\bibitem{CHOST} 
G.~Chen, J.~Han, S.~O, S.~Shan and S.~Tsuchiya,
Finding a spanning Halin subgraph in $3$-connected $\{K_{1,3},P_{5}\}$-free graphs,
preprint.

\bibitem {D}
R.~Diestel,
``Graph Theory'' (4th edition), Graduate Texts in Mathematics \textbf{173},
Springer (2010).

\bibitem{DJG}
D.~Duffus, M.S.~Jacobson and R.J.~Gould,
Forbidden subgraphs and the Hamiltonian theme,
{\it The theory and applications of graphs}, Wiley, New York (1981), 297--316.


\bibitem{FT}
M.~Furuya and S.~Tsuchiya,
Claw-free and $N(2,1,0)$-free graphs are almost net-free,
Graphs Combin. to appear.

\bibitem{H}
R.~Halin,
Studies on minimally $n$-connected graphs,
{\it Combinatorial Mathematics and its Applications}, Academic Press, London (1969), 129--136.

\bibitem{O}
S.~Olariu,
Paw-free graph,
Inform. Process. Lett. \textbf{28} (1988), 53--54.
\end{thebibliography}
\end{document}